\documentclass[12pt,onecolumn,twoside]{article}  

\usepackage{subcaption}
\usepackage{color}
\usepackage{nicefrac}
\usepackage{float}
\usepackage{pbox}
\usepackage{tabularx}
\usepackage[T1]{fontenc}
\usepackage[latin9]{inputenc}
\usepackage{amsmath}
\usepackage{amssymb}
\usepackage{dsfont}
\usepackage{graphicx}
\usepackage{enumerate}
\usepackage{float}
\usepackage{booktabs}
\usepackage{comment}
\usepackage{mathrsfs}
\usepackage{soul}
\usepackage{cite}
\usepackage[margin=1in]{geometry}

\makeatletter

\def\flag={1}
\def\on={1}

\pdfpageheight\paperheight
\pdfpagewidth\paperwidth

\floatstyle{ruled}
\newfloat{algorithm}{tbp}{loa}
\providecommand{\algorithmname}{Algorithm}
\floatname{algorithm}{\protect\algorithmname}
\usepackage[font=small]{caption}

\@ifundefined{showcaptionsetup}{}{%
 \PassOptionsToPackage{caption=false}{subfig}}
\usepackage{subfig}
\makeatother

\def\QED{~\rule[-1pt]{5pt}{5pt}\par\medskip}
\newtheorem{defn}{Definition}
\newtheorem{rem}{Remark}
\newtheorem{lem}{Lemma}

\newtheorem{prop}{Proposition}

\newtheorem{thm}{Theorem}

\newtheorem{coro}{Corollary}
\newtheorem{assumption}{Assumption}
\newenvironment{IEEEproof}{{\bf Proof: }}{ \hfill \QED}
\newtheorem{example}{Example}

\newcommand{\s}{\mathcal{S}}

\newcommand{\Rp}{\mathcal{R}_p}

\newcommand{\R}{\mathbb{R}}

\newcommand{\Hinf}{\mathcal{H}_\infty}
\newcommand{\RHinf}{\mathcal{R}\Hinf}
\newcommand{\A}{\mathscr{A}}
\newcommand{\B}{\mathscr{B}}

\newcommand{\Gs}{\mathcal{A}_*}

\newcommand{\supp}[1]{\text{supp}\left(#1\right)}

\newcommand{\minimize}[1]{\underset{#1}{\text{minimize}}}
\newcommand{\argmin}[1]{\underset{#1}{\text{argmin}}}
\newcommand{\Htwo}{\mathcal{H}_2}

\newcommand{\Htwonorm}[1]{\|#1\|_{\Htwo}}

\newcommand{\Gnorm}[1]{\left\|#1\right\|_{\text{act}}}
\newcommand{\Gdnorm}[1]{\left\|#1\right\|_{\text{act}}^\star}

\newcommand{\apriori}{\emph{a priori}}

\newcommand{\sM}[1]{\sigma_{\max}\left( #1 \right)}
\newcommand{\sm}[1]{\sigma_{\min}\left( #1 \right)}

\newcommand{\drv}{T^{\leq t,v}}
\newcommand{\SNR}[1]{\text{SNR}\left(#1\right)}

\newcommand{\adj}[1]{\left[#1\right]^\dag}
\newcommand{\Noise}{\eta^{\leq t,v}_{\m_*}}

\newcommand{\venkat}[1]{\ifx\flag\on \textbf{\textbf{[{V: #1}]}} \fi}

\newcommand{\cost}[1]{\Psi\left( #1 \right)}
\newcommand{\norm}[1]{\left\| #1 \right\|}

\newcommand{\D}{\mathcal{D}}
\newcommand{\M}{\frak{L}}
\newcommand{\Mtv}{\frak{L}^{\leq t,v}}
\newcommand{\Mtva}{\Mtv_{\Gs}}
\newcommand{\Mtvap}{\Mtv_{\Gs^\perp}}
\newcommand{\atv}{\alpha^{\leq t,v}}
\newcommand{\btv}{\beta^{\leq t,v}}
\newcommand{\gtv}{\gamma^{\leq t,v}}
\newcommand{\Aa}{\mathcal{A}}
\newcommand{\Bb}{\mathcal{B}}

\newcommand{\approxcitations}{\cite{Bon91}}
\newcommand{\ellcitations}{\cite{CRT06,Donoho06}}
\newcommand{\nuccitations}{\cite{fazelThesis,RFP10}}
\newcommand{\distcitations}{\cite{ZDG96, LL11,  S13, SP10, LL12, LD14,M_CDC14_dhinf}}
\newcommand{\Cl}{\mathfrak{L}}
\newcommand{\m}{\mathscr{M}}
\newcommand{\mix}{\tau_{\m_*}}
\newcommand{\wt}{W^{\leq t}}
\newcommand{\tail}{T^{\leq t,v}}
\newcommand{\err}{\wt + \tail}
\newcommand{\note}[1]{{#1}}
\usepackage{tikz}
\usetikzlibrary{arrows}
\usetikzlibrary{shapes}

\newcommand{\edit}[1]{{\textcolor{black}{#1}}}

\renewcommand{\nu}{\delta}

\newcommand{\squeezeup}{\vspace{-2.5mm}}

\newcolumntype{b}{X}
\newcolumntype{s}{>{\hsize=.7\hsize}X}
\newcolumntype{t}{>{\hsize=.3\hsize}X}
\newcolumntype{S}{>{\hsize=.5\hsize}X}

\usepackage{enumitem}
\setlist[itemize]{leftmargin=*}
\setlist[enumerate]{leftmargin=*}

\begin{document}

\title{Regularization for Design}

\author{Nikolai Matni and Venkat Chandrasekaran
\thanks{N. Matni is with the Department of Control and Dynamical Systems, California Institute of Technology, Pasadena, CA.
nmatni@caltech.edu.  V. Chandrasekaran is with the Departments of Computational and Mathematical Sciences and of Electrical Engineering, California Institute of Technology, Pasadena, CA. {venkatc@caltech.edu}.  N.M. was supported in part by the NSF, AFOSR, ARPA-E, and the Institute for Collaborative Biotechnologies through grant W911NF-09-0001 from the U.S. Army Research Office. V.C. was supported in part by NSF Career award CCF-1350590 and by Air Force Office of Scientific Research grant FA9550-14-1-0098. The content does not necessarily reflect the position or the policy of the Government, and no official endorsement should be inferred.  A preliminary version of this work \cite{MC_CDC14} has appeared at the 53rd Annual Conference on Decision and Control in December, 2014.}}
\maketitle
\begin{abstract}
%

When designing controllers for large-scale systems, the architectural aspects of the controller such as the placement of actuators, sensors, and the communication links between them can no longer be taken as given.  The task of designing this architecture is now as important as the design of the control laws themselves.  By interpreting controller synthesis (in a model matching setup) as the solution of a particular linear inverse problem, we view the challenge of obtaining a controller with a desired architecture as one of finding a structured solution to an inverse problem.  Building on this conceptual connection, we formulate and analyze a framework called \emph{Regularization for Design (RFD)}, in which we augment the variational formulations of controller synthesis problems with convex penalty functions that induce a desired controller architecture.  The resulting regularized formulations are convex optimization problems that can be solved efficiently; \edit{these convex programs provide a unified computationally tractable approach for the simultaneous co-design of a structured optimal controller and the actuation, sensing and communication architecture required to implement it.}  Further, these problems are natural control-theoretic analogs of prominent approaches such as the Lasso, the Group Lasso, the Elastic Net, and others that are employed in statistical modeling.  In analogy to that literature, we show that our approach identifies optimally structured controllers under a suitable condition on a ``signal-to-noise'' type ratio.

\end{abstract}

\section{Introduction}
\label{sec:intro}

As we move into the era of large-scale systems such as the smart-grid, software defined networking and automated highways, the design of control systems is becoming increasingly more challenging.  Designing the controller \emph{architecture} -- placing sensors and actuators as well as the communication links between them -- is now as important as the traditional design of the control laws themselves.

A conceptually useful viewpoint in the design of controller architectures is to consider complicated systems as being composed of multiple simpler \emph{atomic} subsystems.  For example, if the task is to design the actuation architecture of a controller, a natural atomic element is a controller with a single actuator -- it is then clear that a general architecture can be built out of such atoms.  In general, controllers with a dense actuation, sensing and communication architecture (i.e., systems that consist of many atomic subsystems) achieve better closed loop performance in comparison with those with sparse architectures (i.e., systems composed of a small number of atomic subsystems).  However, as these architectural resources translate into actual hardware requirements, it is desirable from both a maintenance and a cost perspective that we minimize the total number of atomic elements used.  Hence, the problem of controller architecture/control law co-design is one of jointly optimizing an appropriately defined structural measure of the controller and its closed loop performance by trading off between these two competing metrics in a principled manner.  In other words, we seek an approximation of a given optimal controller by one that utilizes fewer atomic elements without a significant loss in performance.  This goal has parallels with the approximation theory literature in which one seeks approximations of complicated functions as combinations of elements from a simpler class of functions such as the Fourier basis or a wavelet basis \approxcitations.

In an appropriate parameterization, pure controller synthesis methods in a model matching framework can be interpreted as techniques for solving a particular linear inverse problem in which one is given an open loop response of a system $Y$ and a linear map $\Cl$ from the controller to the closed loop response, and one seeks a controller $U$ such that \edit{$Y - \Cl(U) \approx 0$} (as measured in a suitable performance metric).  From this perspective, our objective in joint controller architecture/control law co-design is to obtain \emph{structured} solutions to the linear inverse problem underlying controller synthesis.  Such \emph{structured linear inverse problems} (SLIP) are of interest in diverse applications across applied mathematics -- for instance, computing sparse solutions to linear inverse problems or computing low-rank solutions to systems of linear matrix equations arise prominently in many contexts in signal processing and in statistics \cite{CRT06,Donoho06,fazelThesis,RFP10}.  

In these problem domains, minimizing the $\ell_1$ norm subject to constraints described by the specified equations is useful for obtaining sparse solutions \ellcitations, and similarly, nuclear norm minimization is useful for obtaining low-rank solutions to linear matrix equations \nuccitations.  These ideas were extended in \cite{CRPW12}, where the authors describe a generic convex programming approach -- based on minimizing an appropriate \emph{atomic norm} \approxcitations  -- for inducing a desired type of structure in solutions to linear inverse problems.  Motivated by these developments, our approach to the problem of joint architecture/control law co-design is to augment variational formulations of controller synthesis methods with suitable convex regularization functions.  The role of these regularizers is to penalize controllers with more complex architectures in favor of those with less complex ones, thus inducing controllers with a simpler architecture.  We call this framework \emph{Regularization for Design} (RFD).

\textbf{Related Work:} Regularization techniques based on $\ell_1$ norms and, more generally, atomic norms have already been employed extensively in system identification, e.g., to identify systems of small Hankel order (cf. \cite{SBTR12,Fazel,MR_CDC14}), and in linear regression based methods \cite{LjungNewOld}.  Although the resulting solutions yield structured systems, they typically do not have a direct interpretation in terms of the architecture of a control system (i.e., actuators, sensors and the communication links between them).  The use of regularization explicitly for the purpose of designing the architecture of a controller can also be found in the literature.  Examples include the use of $\ell_1$ regularization to design sparse structures in $\mathcal{H}_2$ static state feedback gains \cite{LFJ13}, treatment therapies \cite{JRM14}, and synchronization topologies \cite{FLJ13}; the use of group norm penalties to design actuation/sensing schemes \cite{MC_CDC14,DJL14}; and the use of a specialized atomic norm to design communication delay constraints that are well-suited to $\mathcal{H}_2$ distributed optimal control \cite{M_CDC13_codesign, M_TCNS14}.  Although these methods provide an algorithmic approach for designing controllers with a desired architecture in certain specialized settings, they do not enjoy the same theoretical support that regularization techniques for structured inverse problems enjoy in other settings \cite{Bon91,CRT06,Donoho06,fazelThesis,RFP10}.  \edit{Methods based on greedy algorithms have also been developed to identify minimal architectures that guarantee the structural observability and controllability of a system \cite{PKP15}.}


\textbf{Our Contributions:} 
\edit{This paper presents novel computational and theoretical contributions to the area of optimal controller architecture/control law co-design.  From a computational perspective, we propose a general RFD framework that is applicable in a much broader range of settings than the previous approaches mentioned above.  We restrict ourselves to problems for which the linear optimal \note{structured} controller is specified as the solution to a convex optimization problem \cite{RL06,RCL10,LL11_QI}.  As a result, RFD optimization problems with convex regularization functions for inducing a desired architecture are convex programs.   Specifically,
(i) we provide a catalog of atomic norms useful for control architecture design.  In particular, in addition to known penalties for actuator, sensor and communication design, we provide novel penalties for \emph{simultaneous} actuator, sensor and/or communication design;
(ii) we describe a unifying framework for RFD that encompasses state and output feedback problems in centralized and distributed settings, and in which any subset of actuation, sensing, and/or communication architectures are co-designed; and
(iii) we present a two-step algorithm that first identifies the controller architecture via a finite-dimensional convex RFD optimization problem, and then solves for the potentially infinite dimensional linear optimal controller restricted to the designed architecture using methods from classical and distributed optimal control \distcitations.}

\edit{To provide theoretical support for our proposed computational framework, we make explicit links between RFD optimization problems and the use of convex optimization based approaches for structured inference problems.  We build on these links to analyze the properties of the structured controllers generated by RFD synthesis methods, which leads to conditions under which RFD methods succeed in identifying optimally-structured controllers.  
Our analysis and results are natural
control-theoretic analogs of similar results in the statistical
inference literature.
Specifically, (i) 
we show that finite-horizon finite-order convex approximations of an RFD optimization problem can recover the structure of an underlying infinite dimensional optimal controller; (ii) we define control-theoretic analogs of identifiability conditions and signal-to-noise ratios (SNRs), and we provide sufficient conditions based on these for a controller architecture to be identified by RFD. In particular, we show that controllers that maximize this SNR-like quantity are more easily recovered via RFD than those that do not, and (iii) we provide a concrete example of a system satisfying the above identifiability and SNR conditions.  As far as we are aware, this is the first example in the literature of a system for which it is shown that convex optimization provably recovers the actuation architecture of an underlying optimally structured controller.
}

%
%
%

\textbf{Paper Organization}: \edit{In  \S\ref{sec:prelims}, we define notation and discuss the relevant concepts from operator theory.  This paper is then organized in a modular fashion: \S\ref{sec:RFD}-\ref{sec:rfd_proc} focus on the computational aspects of controller architecture design, whereas \S\ref{sec:recovery} and \S\ref{sec:case} focus on conditions for optimal architecture recovery.  Specifically, in \S\ref{sec:RFD}, we introduce the RFD framework as a natural blend of controller synthesis methods and regularization techniques.  In \S\ref{sec:RFD_algos}, we focus on RFD problems with an $\Htwo$ performance metric and an atomic norm penalty; we present a catalog of atomic norms that are useful for controller architecture design, and make connections to the structured inference literature.  The computational component of the paper concludes in \S\ref{sec:rfd_proc}, in which we formally describe the two-step RFD procedure, and we apply the RFD framework to a simultaneous actuator, sensor and communication design problem.  In \S\ref{sec:recovery}, we shift our focus to analyzing the theoretical properties of the RFD procedure:  we make connections between structured controller design and structured inference problems by framing both tasks as finding structured solutions to linear inverse problems, and we leverage these connections to describe interpretable sufficient conditions for the success of a finite-dimensional RFD optimization problem.  In  \S\ref{sec:case}, we provide a case study to further illustrate the applicability of these results. 
}

\section{Preliminaries \& Notation}
\label{sec:prelims}
\label{sec:notes}
We use standard definitions of the Hardy spaces $\Htwo$ and $\Hinf$.  As is standard, we denote the restrictions of $\Htwo$ and $\Hinf$ to the space of real proper transfer matrices $\Rp$ by $\mathcal{R}\Htwo$ and $\RHinf$.  As we work in discrete time, the two spaces are equal, and as a matter of convention  we refer to this space as $\RHinf$.   We refer the reader to \cite{ZDG96} for a review of this standard material.
\edit{
Let $\RHinf^{\leq t}$ denote the subspace of $\RHinf$ composed of finite impulse response (FIR) transfer matrices of length $t$, i.e., $\RHinf^{\leq t} := \{ G \in \RHinf \, | \, G = \sum_{i=0}^t \frac{1}{z^i} G^{(i)} \}$.  We denote the projection of an element $G \in \Rp$ onto the subspace $\RHinf^{\leq t}$ by $G^{\leq t}$.
}
Unless required for the discussion, we do not explicitly denote dimensions and we assume that all vectors, operators and spaces are of compatible dimension throughout.

\edit{We denote elements of $\Rp$ with upper case Latin letters, and temporal indices and horizons by lower case Latin letters.  Linear maps from $\Rp$ to $\Rp$ are denoted by upper case Fraktur letters such as $\M$.  For such a linear map, we denote the $i^{\textrm{th}}$ impulse response element of $\M$ by $\M^{(i)}$.  We further use $\M^{\leq t}$ to denote the restriction of the range of $\M$ to $\RHinf^{\leq t}$, and $\M^{\leq t,v}$ to denote the restriction of $\M^{\leq t}$ to the domain $\RHinf^{\leq v}$.  Thus $\M^{\leq t,v}$ is a map from $\RHinf^{\leq v}$ to $\RHinf^{\leq t}$.  In particular, if $\M$ is represented as a semi-infinite lower block triangular matrix, then $\M^{\leq t,v}$ corresponds to the $t$ by $v$ block row by block column sub matrix $(\M)_{ij}, \, i=1,\dots,t,\, j = 1,\dots,v$.
}

Sets are denoted by upper case script letters, such as $\mathscr{S}$, whereas subspaces of an inner product space are denoted by upper case calligraphic letters, such as $\s$. The restriction of a linear map $\M$ to a subspace $\s\in\RHinf$ is denoted by $\M_\s$; \edit{ similarly, the projection of an operator $G \in \Rp$ onto a subspace $\Aa\subset\Rp$ is denoted by $G_\Aa$. } We denote the adjoint of a linear map $\M$ by $\M^\dag$.  The most complicated expression that we use is of the form $\left[\M^{\leq t,v}_{\s}\right]^\dag$:  this denotes the adjoint of the map $\M^{\leq t}$ restricted to $\RHinf^{\leq v} \cap \s$.

We denote the $n$-dimensional identity matrix and down-shift matrices by $I_n$ and $Z_n$, respectively.  In particular, $Z_n$ is a matrix with all ones along its first sub-diagonal and zero elsewhere.  We use $e_i$ to denote a standard basis element in $\R^n$, and $E_{ij}$ to denote the matrix with $(i,j)^\text{th}$ element set to 1 and all others set to 0.

\section{RFD as Structured Approximation}
\label{sec:RFD}
\note{Under standard assumptions  \cite{ZDG96},} traditional controller synthesis methods within the framework of \emph{model matching} can be framed as linear inverse problems of the form
\begin{equation}
\begin{array}{rl}
\minimize{U\in \RHinf} & \cost{U;Y,\Cl}\\
\end{array}
\label{eq:trad_opt}
\end{equation}
where $Y$ is the open loop response of the system, $U$ is the Youla parameter, $\Cl$ is \note{a suitably defined} linear map from the Youla parameter $U$ to the closed loop response, and $\cost{\cdot;Y,\Cl}$ is a performance metric \edit{that measures the size of the closed loop response (i.e., the size of the deviation between $Y$ and $\Cl(U)$), such as the $\Htwo$ or $\Hinf$ norm.}  
We make this connection clear in the following subsection and we \edit{also recall how} to incorporate \emph{quadratically invariant} \cite{RL06} distributed constraints on the controller into this framework.  
\edit{
\begin{rem}
We use this non-standard notation to facilitate comparisons with the structured inference literature.  In particular, this notation emphasizes that the optimal linear controller synthesis task can be viewed as one of solving a linear inverse problem with ``data'' specified by the open loop response $Y$ of the system and the map $\Cl$ from the Youla parameter to the closed loop response.
\end{rem}}

\subsection{Convex Model Matching}
\label{sec:model_match}

\begin{figure}[h!]
\centering
\includegraphics[width = .35\textwidth]{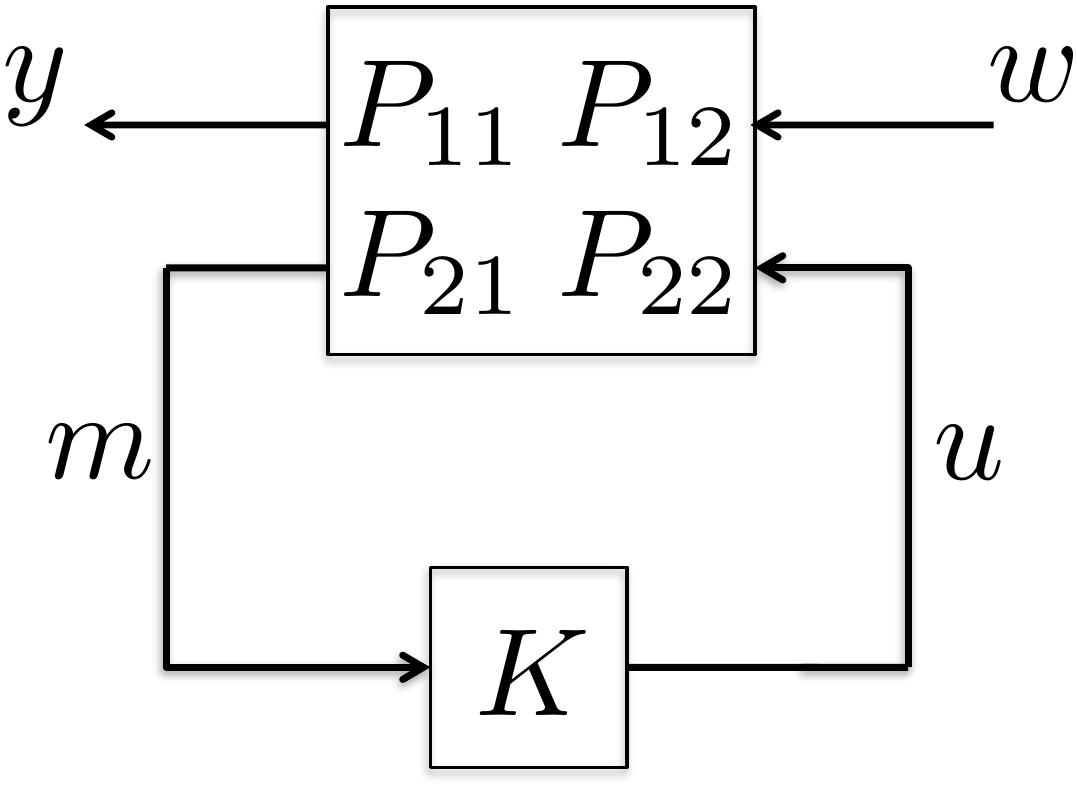}
\caption{A diagram of the generalized plant defined in \eqref{eq:genplant}.}
\label{fig:genplant}
\end{figure}
In order to discuss a broad range of model matching problems, we introduce the \emph{generalized plant}, a standard tool in robust and optimal control \cite{ZDG96}.  In particular, consider the system described by
\begin{equation}
P = \begin{bmatrix}
P_{11} & P_{12} \\
P_{21} & P_{22}
\end{bmatrix}
= \left[\begin{array}{c|cc}
A & B_{1} & B_{2} \\ \hline
C_{1} & 0 & D_{12} \\
C_{2} & D_{21} & 0
\end{array}\right]
\label{eq:genplant}
\end{equation}
where $P_{ij} = C_{i}(zI-A)^{-1}B_{j} + D_{ij}$.  As illustrated in Figure \ref{fig:genplant}, this system describes the four transfer matrices from the disturbance and control inputs $w$ and $u$, respectively, to the controlled and measured outputs $y$ and $m$, respectively.  We make the (slightly modified) standard orthogonality assumptions that 
\begin{equation}
\begin{array}{rcl}
D_{12}\begin{bmatrix} C_1^\top & D_{12}^\top \end{bmatrix} &=&  \begin{bmatrix} 0 & \rho_u I \end{bmatrix}\\
 D_{21}^\top\begin{bmatrix} B_1 & D_{21} \end{bmatrix} &=&  \begin{bmatrix} 0 & \rho_w I \end{bmatrix}
\end{array}
\label{eq:ortho}
\end{equation}
 for some $\rho_u, \rho_w \geq 0$. \note{ At times we separate the state component of the open and closed loop responses from the components of these transfer matrices that measure control effort. To that end, we define the \emph{state component} of an element $X$ to be the projection of $X$ onto the range of $C_1$.}
 
 \edit{Letting $u(z) = K(z)m(z)$ for a causal linear controller $K\in\Rp$, the closed loop map from the disturbance $w$ to the controlled output $y$ is given by the linear fractional transform $P_{11}-P_{12}K(I-P_{22}K)^{-1}P_{21}.$  A typical optimal control problem in this framework is then formulated as 
\begin{equation}
\begin{array}{rl}
\minimize{K\in\Rp} & \left\| P_{11}-P_{12}K(I-P_{22}K)^{-1}P_{21} \right\| \\
\text{s.t.} & K(I-P_{22}K)^{-1} \in \RHinf
\end{array}
\label{eq:OC}
\end{equation}
where $\|\cdot\|$ is a suitable norm, and the constraint ensures internal stability of the closed loop system \cite{ZDG96}.  Notice however that the optimal control problem \eqref{eq:OC} is non-convex as stated.}


\edit{When the open-loop plant is stable (we remark on the unstable case at the end of the subsection), a standard and general approach to solving the optimal control problem \eqref{eq:OC} is to convert it into a \emph{model matching problem} through the \emph{Youla} change of variables $U:= K(I-P_{22}K)^{-1}$; the optimal controller $K$ can then be recovered via $K = (I+UP_{22})^{-1}U$.  The resulting convex optimization problem is then given by
\begin{equation}
\begin{array}{rl}
\minimize{U\in\RHinf} & \left\| P_{11}-P_{12}UP_{21} \right\|, \\
\end{array}
\label{eq:MM}
\end{equation}
which is of the form of the linear inverse problem \eqref{eq:trad_opt} if we take $Y:= P_{11}$, $\Cl= P_{12}\otimes P_{21}$ (where $\left(P_{12}\otimes P_{21}\right)(U) := P_{12} U P_{21}$) and $\cost{U; Y,\Cl} := \left\| Y- \Cl(U)\right\|$.
}



We also often want to impose a \emph{distributed constraint}  on the controller $K$ by requiring $K$ to lie in some subspace $\s$, which specifies information exchange constraints between the sensors and actuators of the controller.  It is known that a necessary and sufficient condition for such a distributed constraint to be invariant under the Youla change of variables is that it be \emph{quadratically invariant} with respect to $P_{22}$ \cite{RL06,RCL10,LL11_QI}.

\begin{defn}[Quadratic Invariance]
A subspace $\s$ is \emph{quadratically invariant} (QI) with respect to $P_{22}$ if $KP_{22}K \in \s \ \forall K \in \s$.
\label{def:QI}
\end{defn}

In particular, when a subspace $\s$ is QI with respect to $P_{22}$, we have that $K \in \s$ if and only if $K(I-P_{22}K)^{-1} \in \s$, allowing us to convert the general optimal control problem \eqref{eq:trad_opt} with the additional constraint that $K\in\s$ to the following convex optimization problem

\begin{equation}
\begin{array}{rl}
\minimize{U\in\RHinf} & \cost{U; Y, \Cl} \
\text{s.t. } U \in \s.
\end{array}
\label{eq:gen_mm}
\end{equation}
\note{
This optimization problem is again precisely of the form of the linear inverse problem \eqref{eq:trad_opt} save for the addition of the subspace constraint $U \in \s$.}
\edit{This framework is fairly general in that it allows for a unified treatment of all structured optimal control problems in which the linear optimal structured controller can be computed via convex optimization \cite{RL06,LL11_QI}.  These include state and output feedback problems in either centralized or QI distributed settings.  \edit{Further, if the optimal control problem is centralized with respect to the $\Htwo$, $\Hinf$ or $\mathcal{L}_1$ metrics, or is QI distributed with a finite horizon $\Htwo$ cost, the linear optimal control is globally optimal \cite{ZDG96,HC72,DahlehL1,LL11_QI}.}}
\edit{
\begin{rem}[Extension to Unstable Plants]
The above discussion extends to unstable plants through the use of an appropriate \emph{structure preserving} Youla-Kucera parameterization \edit{built around \emph{arbitrary coprime factorizations}}, which are always available.   See \cite{SM12,MMRY12,LD14} for examples of such parameterizations, and \cite{M_TCNS14} for an example of using such a parameterization with a structure inducing penalty.
\label{rem:unstable}
\end{rem}}

%
\vspace{-2mm}

\subsection{Architecture Design through Structured Solutions}
\label{sec:RFDQI}

We seek a modification of the optimal controller synthesis procedure to design the controller's architecture.  We reiterate that by
the {architecture} of a controller, we mean the actuators, sensors, and communication links between them.  In particular, we view the controller $K$ as a map from all potential sensors to all potential actuators, using all potential communication links between these actuators and sensors.  The architectural design task is that of selecting which actuators, sensors and communication links need to be used to achieve a certain performance level.  This task is naturally viewed as one of finding a {structured approximation} of the optimal controller \eqref{eq:OC} that utilizes all of the available architectural choices.  

The components of the controller architecture being designed determine the type of structured approximation that we attempt to identify. In particular, each nonzero row of $K(z)$ corresponds to an actuator used by the controller, and likewise, each nonzero column corresponds to a sensor employed by the controller.  Further sparsity patterns present within rows/columns of the power series elements $K^{(t)}$ of $K(z)$ can be interpreted as information exchange constraints imposed by an underlying communication network between the sensors and actuators.  With these observations, it is clear that specific sparsity patterns in $K$ have direct interpretations in terms of the architectural components of the controller: nonzero rows correspond to actuators, nonzero columns correspond to sensors, and additional sparsity structure corresponds to communication constraints.   

\note{Although we seek seek to identify a suitably structured controller $K$, for the computational reasons described in \S\ref{sec:model_match} it is preferable to solve a problem in terms of the Youla parameter $U$ as this parameterization leads to the convex optimization problem \eqref{eq:gen_mm}.  
Therefore, in the following definition RFD problems are defined as a regularized version of the model matching problem \eqref{eq:gen_mm} with a penalty function added to the objective to induce suitable structure, rather than as a modification of the controller synthesis problem \eqref{eq:OC}.   In the sequel, we justify that for architectural design problems of interest, the structure underlying the controller $K$ is equivalent to the structure underlying the Youla parameter $U$;   to that end, we show that in the case of actuator, sensor, and/or QI communication topology design,\footnote{\edit{We restrict ourselves to communications delays that satisfy the triangle inequality defined in \cite{RCL10}. This assumption implies that information exchanged between sensors and actuators is transmitted along shortest delay paths in the underlying communication graph.}} the structure of the Youla parameter $U$ corresponds to the structure of the controller $K$.}


%

\begin{defn}
Let $U, \, Y\in\RHinf$, and $\Cl:\RHinf \to \RHinf$ be of compatible dimension.  The optimization
\begin{equation}
\begin{array}{rl}
\minimize{U\in \RHinf} & \cost{U; Y, \Cl}+2\lambda\Omega\left(U\right) \
\text{s.t. } U \in \s
\end{array}
\label{eq:rfd_opt}
\end{equation}
is called a \emph{RFD optimization problem} with cost function $\cost{\cdot;Y,\Cl}$ and penalty $\Omega\left(\cdot\right)$.
\end{defn}

\begin{rem}[Static or Dynamic]
With the exception of the centralized state-feedback setting, it is known that the optimal linear controller is dynamic \cite{ZDG96}, and therefore restricting our analysis to dynamic controllers is natural.  In the centralized setting (given the equivalence between static and dynamic state-feedback), once an appropriate actuation architecture\footnote{For a system to be centralized, there can be no communication constraints, and for it to be state-feedback, there can be no sensing constraints, leaving only actuation architecture design as a possible RFD task.} is identified, traditional methods can then be used to solve for a static state-feedback controller restricted to that architecture.  Further as we show in \S\ref{sec:recovery}, the dynamic controller synthesis approach is amenable to analysis that guarantees {optimal structure recovery}.
\end{rem}

We discuss natural costs $\cost{\cdot;Y,\Cl}$ and penalties $\Omega\left(\cdot\right)$ in \S\ref{sec:RFD_algos}, and we focus now on justifying why we can perform the structural design on the Youla parameter $U$ rather than the controller $K$. 
\edit{
\paragraph{Actuator/Sensor Design}  \note{Recall that the actuators (sensors) that a controller $K$ uses are identified by the nonzero rows (columns) of $K$: the actuator (sensor) design problem therefore corresponds to finding a controller $K$ that achieves a good closed loop response and that is sparse row-wise (column-wise).  This corresponds exactly to finding a row (column) sparse solution $U$ to the RFD optimization problem \eqref{eq:rfd_opt}.  This is true because any subspace $\D$ that is defined solely in terms of row (column) sparsity is QI with respect to any $P_{22}$.   In particular, it is easily verified that if $K \in \D$, then right (left) multiplication leaves $\D$ invariant, i.e., $KX \in \D$  ($XK\in \D$) for all compatible $X$. It then follows from Definition \ref{def:QI} that $\D$ is QI with respect to any plant $P_{22}$.  We can extend this analysis to incorporate additional QI distributed constraints $\s$ by leveraging the results in \cite{RCL10}: in particular, if $\s$ is QI with respect to $P_{22}$, then so is $\D \cap \s$.}}\footnote{In particular, since removing actuators does not change the communication delays between the remaining actuators and sensors, if the delay based conditions in \cite{RCL10}  hold when all actuators (sensors) are present they also hold with any subset of them being present.}


\paragraph{Joint Actuator and Sensor Design} \note{By virtue of the previous discussion joint actuator and sensor design corresponds to finding a controller $K$ that is simultaneously sparse row-wise and column-wise.  It follows immediately from the previous discussion that this corresponds exactly to finding a simultaneously row and column sparse solution $U$ to the RFD optimization problem \eqref{eq:rfd_opt}.  In particular, any subspace $\D$ defined solely in terms of row and column sparsity is QI with respect to any plant $P_{22}$, we can incorporate additional QI distributed constraints $\s$ by leveraging the results in \cite{RCL10}.}

\paragraph{Communication Design}
\note{In an analogous manner to the above, one can also associate subspaces to structures corresponding to suitable information exchange constraints that a distributed controller must satisfy.
Recall in particular that the information exchange constraints between the sensors and actuators of a controller $K$ are identified by the sparsity structure found within the nonzero rows and columns of $K$.  In \cite{M_TCNS14}, the first author showed that a specific type of sparsity structure in $K$ corresponds exactly to sensors and actuators exchanging information according to an underlying communication graph.
In particular, given a communication graph between sensors and actuators with adjacency matrix $\Gamma$, a distributed controller $K$ can be implemented using the graph defined by $\Gamma$ if the power series elements $K^{(t)}$ of the controller satisfy $\supp{K^{(t)}}\subseteq \supp{\Gamma^{t-1}}$.\footnote{We assume that $K$ is square for simplicity; cf. \cite{M_TCNS14} for the general case.}  The interpretation of the support nesting condition is that the delay from sensor $j$ to actuator $i$ is given by the length of the shortest path from node $j$ to node $i$ in the graph defined by $\Gamma$.  This support nesting condition thus defines the distributed subspace constraint $\s$ in which $K$ must lie -- based on the discussion in \S\ref{sec:model_match}, one can pose the distributed controller synthesis problem as a distributed model matching problem \eqref{eq:gen_mm} if and only if $\s$ is QI with respect to $P_{22}$.  In light of this, we consider the communication design task proposed in \cite{M_TCNS14}: given an initial graph with adjacency matrix $\Gamma_{\text{QI}}$ that induces a QI distributed subspace constraint $\s$, what minimal set of additional edges should be added to the graph to achieve a desired performance level.\footnote{It is shown in \cite{M_TCNS14} that under mild assumptions on the plant $P_{22}$, the propagation delays of $P_{22}$ can be used to define an adjacency matrix $\Gamma_{\text{QI}}$ that induces a distributed subspace constraint that is QI with respect to $P_{22}$.}  It is additionally shown in \cite{M_TCNS14} that any communication graph constructed in this manner results in a subspace constraint $\s$ that is QI with respect to $P_{22}$.  Therefore the structure imposed on the controller $K$ by an underlying QI communication graph corresponds exactly to the structure imposed on the Youla parameter $U$.
 }
 

\paragraph{Joint Communication, Actuator and/or Sensor Design}
\note{By virtue of the previous discussion and the results of \cite{RCL10}, combining QI communication design with actuator and/or sensor design still leads to the underlying structure of the controller $K$ corresponding to the underlying structure of the Youla parameter $U$.}

Thus for architecture design problems the RFD task can be performed via a model matching problem.

\begin{example}
Suppose that different RFD optimization problems are solved for a system with three possible actuators and three possible sensors, resulting in the various sparsity patterns in $U(z)$ shown on the far right of Figure \ref{fig:patterns}.  It is easily seen by inspection that the resulting sparsity patterns are QI. In particular Figures \ref{fig:patterns}a) through \ref{fig:patterns}c) correspond to centralized RFD optimization problems (this can be seen from the full matrices in the center of the left hand side), and 2d) to a RFD optimization problem subject to lower triangular constraints, a special case of a nested information constraint.\footnote{Nested information constraints are a well studied class of QI distributed constraints, cf. \cite{LL11, S13, SP10} for examples.}

\begin{figure}[h]
\centering
\includegraphics[width = .5\textwidth]{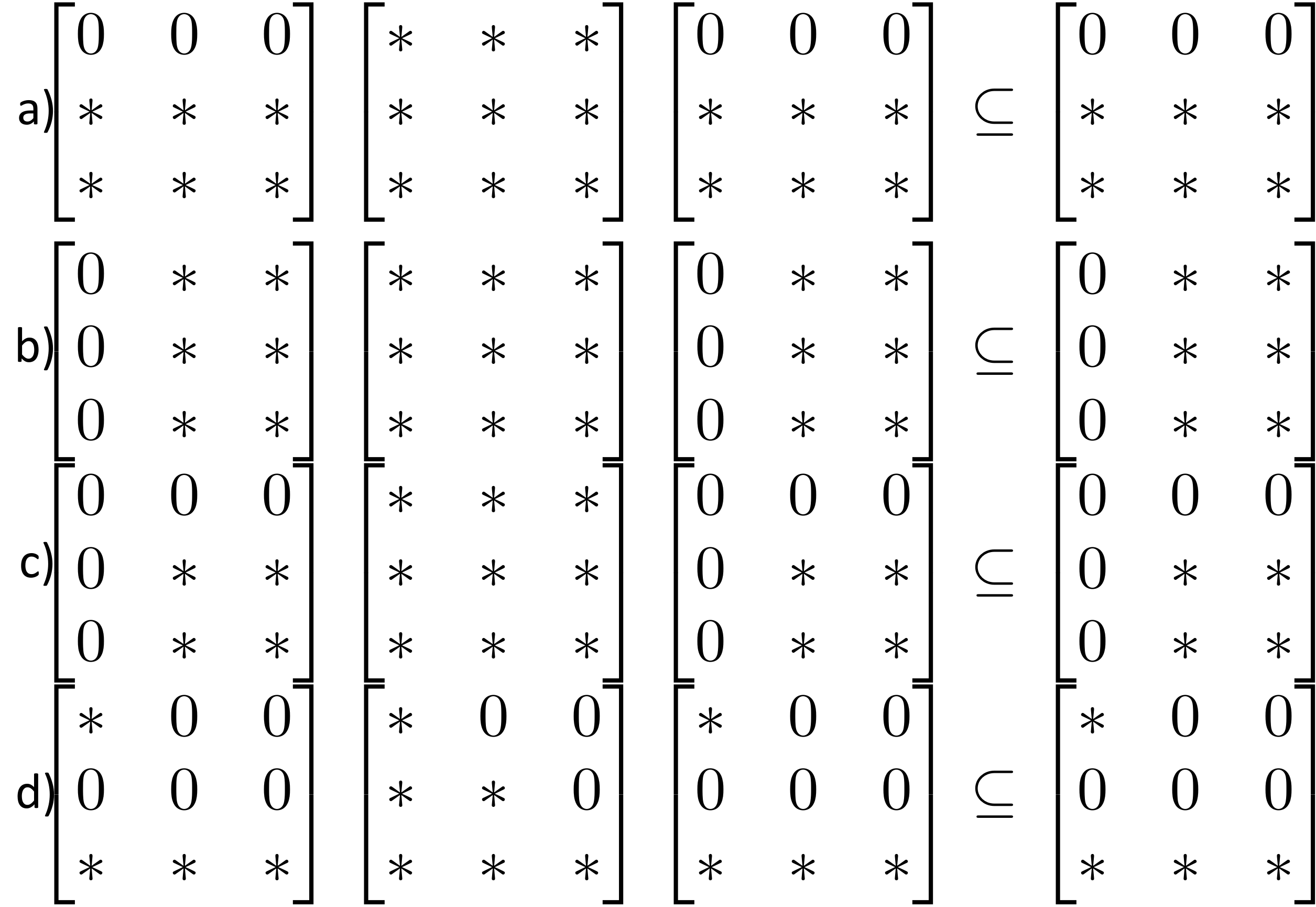}
\caption{Examples of QI sparsity patterns generated via a) actuator, b) sensor, and c) actuator/sensor RFD procedures without any distributed constraints, and d) actuator RFD subject to nested information constraints.}
\label{fig:patterns}
\end{figure}
\end{example}
\squeezeup \squeezeup \squeezeup

\section{RFD Cost Functions and Regularizers}
\label{sec:RFD_algos}
In this section we examine convex formulations of the RFD optimization problem \eqref{eq:rfd_opt} by restricting our attention to convex cost functions $\cost{\cdot; Y, \Cl}$ and {convex} penalty functions $\Omega\left(\cdot\right)$.
\squeezeup
\subsection{Convex Cost Functions}
\label{sec:cost_fcns}
Any suitable convex cost function $\cost{\cdot; Y, \Cl}$ can be used in \eqref{eq:rfd_opt}: traditional examples from robust and optimal control include the $\Htwo$, $\Hinf$ \cite{ZDG96} and $\mathcal{L}_1$ norms \cite{DahlehL1}.  We focus on the $\Htwo$ norm as a performance metric because it allows us to make direct connections between the RFD optimization problem \eqref{eq:rfd_opt} and well-established methods employed in statistical inference such as ordinary least squares, Ridge Regression \cite{Ridge}, Group Lasso \cite{GroupLasso} and Group Elastic Net \cite{ElasticNet}.

We begin by introducing a specialized form of the model matching problem \eqref{eq:gen_mm}, and show how state-feedback problems with $\Htwo$ performance metrics can be put into this form.

\begin{defn}
Let $U, \, Y\in\RHinf$, and $\M:\RHinf \to \RHinf$ be of compatible dimension.  The optimization problem
\begin{equation}
\begin{array}{rl}
\minimize{U\in \RHinf} & \Htwonorm{Y - \M(U)}^2+\rho_u\Htwonorm{U}^2 \
\text{s.t. }  U \in \s
\end{array}
\label{eq:the_opt}
\end{equation}
is the $\Htwo$ \emph{optimal control problem} for a suitable control penalty weight $\rho_u$ and a distributed constraint $\s$.
\end{defn}

\edit{In this definition, $Y$ is the state component of the open loop response, and $\M$ is the map from the Youla parameter to state component of the closed loop response. }
Explicitly separating the cost of the state component $\Htwonorm{Y - \M(U)}^2$ of the closed loop response from the control cost $\rho_u\Htwonorm{U}^2$ of the closed loop response allows us to connect the $\Htwo$ {RFD optimization} to several well-established methods in the inference literature.  
Before elaborating on some of these connections, we provide two examples of standard control problems that can be put into this form.

\begin{example}[Basic LQR]
Consider the basic LQR problem given by
\begin{equation}
\begin{array}{rl}
\minimize{u\in \ell_2} & \sum_{t=0}^\infty \|Cx_t\|_{\ell_2}^2 + \|Du_t\|_{\ell_2}^2 \\
\text{s.t.} & x_{t+1} = Ax_t + Bu_t,\ x_0 = \xi,
\end{array}
\end{equation}
and assume that $D^\top D = \rho_u I$, for some $\rho_u\geq0$.
Define $\rho = \rho_u$, $X^{(t)} = CA^{t}\xi$ for $t\geq 0$, $U^{(t)} = u_t$, and $\M(U) = -H\ast U$, where $H\in \Rp$ with $H^{(0)}=0$, and $H^{(t)} = CA^{t-1}B$ for $t\geq 1$.  We can then rewrite the basic LQR problem in the form of optimization problem \eqref{eq:the_opt} (with no distributed constraint $\s$).
\label{ex:bLQR}
\end{example}
\squeezeup \squeezeup
\edit{
\begin{example}[$\Htwo$ State Feedback]
Assume either that the generalized plant \eqref{eq:genplant} is open-loop stable or that the control problem is over a finite horizon.  Let $C_2 = I$ and $D_{21} =0$ in the generalized plant \eqref{eq:genplant} such that the problem is one of synthesizing an optimal state-feedback controller, and for clarity of exposition, assume that $B_1$ is invertible.\footnote{\edit{The assumption that $B_1$ is invertible simply implies that no component of the state is deterministic, and can be relaxed at the expense of more complicated formulas.}}   Define the Youla parameterization for the controller synthesis problem \eqref{eq:OC} as follows \cite{ZDG96}:
\begin{equation*}
\begin{array}{l}
\tilde{P}_{12} = \frac{1}{z}P_{12}, \
\tilde{P}_{21} = AP_{21} + B_1, \
\tilde{U}=K(I-P_{22}K)^{-1}\tilde{P}_{21},
\end{array}
\end{equation*}
with all other parameters remaining the same.  Under this parameterization, the optimal control problem \eqref{eq:OC} (with additional QI distributed constraint $\s$) with performance metric $\Htwonorm{\cdot}^2$ can be written as
\begin{equation}
\begin{array}{l}
\minimize{\tilde{U}\in\RHinf} \Htwonorm{P_{11} - \M\left(\tilde{U}\right)}^2 + \rho_u\Htwonorm{\tilde{U}}^2 \
\text{s.t. }  \tilde{U}\tilde{P}_{21}^{-1} \in \s,
\end{array}
\label{eq:sf_mm}
\end{equation}
The optimal controller $K$ is then recovered from the solution to \eqref{eq:sf_mm} as $K=(I+\tilde{U}\tilde{P}_{21}^{-1}P_{22})^{-1}\tilde{U}\tilde{P}_{21}^{-1}$.  The state-feedback assumption and the choice of Youla parameterization ensure that $\tilde{P}_{21}$ is invertible in $\RHinf$ and that $K \in \s$.  
\begin{rem}
A dual argument applies to $\Htwo$ filter design by considering the ``full-control'' setting, cf. \cite{ZDG96} for more details.
\end{rem}
\label{ex:actuators}
\end{example}}

\edit{As illustrated by Example \ref{ex:actuators}, the $\Htwo$ optimal control problem \eqref{eq:the_opt} is simply a more general way of writing the $\Htwo $ state feedback model matching problem -- in Remark \ref{rem:output_feedback} we show how $\Htwo$ output feedback model matching problems can also be put in a similar form.  Writing the $\Htwo$ problem as a linear inverse problem with a least squares like state cost and an explicitly separated control cost already allows us to make connections to classical techniques from statistical inference.  These connections (along with others we make later in this
section) are summarized in Table \ref{tab:dict}.}  In order to keep the discussion as streamlined as possible, we make these connections in the context of the Basic LQR problem presented in Example \ref{ex:bLQR}.

\subsubsection{$\rho_u=0$}  In an inferential context, this is simply ordinary least squares, and is commonly used when $U_*$ is not known \emph{a priori} to have any structure.  Further, the resulting estimate of $U_*$ is unbiased, but often suffers from high error variance.  Moving now to a control context, It is easy to see that this setting corresponds to ``cheap control'' LQR, in which there is no cost on $u_t$ -- under suitable \edit{controllability and observability} assumptions, the resulting state trajectory is deadbeat, but the optimal control law is not necessarily unique.

\subsubsection{$\rho_u >0$}  \note{This corresponds to Ridge Regression or Tikhonov Regularization \cite{Ridge}.  In an inferential context, this regularizer has the effect of \emph{shrinking} estimates towards 0 -- this introduces bias into the estimator, but reduces its variance, and is often a favorable tradeoff from a statistical perspective.  From a linear algebraic perspective, this is a commonly used technique to improve the numerical conditioning of an inverse problem.  Once again, the interpretation in RFD is clear: this corresponds to standard LQR control with $R = \rho I$; the parameter $\rho$ allows for a tradeoff between control effort and state deviation.  The optimal control action is then {unique} and the resulting state trajectory is generally not deadbeat.}
\squeezeup
\subsection{The $\Htwo$ RFD Optimization Problem with an Atomic Norm Penalty}
\label{sec:atoms}
Recall that our strategy for designing controller architectures is to augment the traditional model matching problem \eqref{eq:gen_mm} with a structure inducing penalty, resulting in the RFD optimization problem \eqref{eq:rfd_opt}.  In light of the discussion of the previous subsection, we further specialize the RFD optimization problem to have an $\Htwo$ performance metric and an atomic norm penalty $\norm{\cdot}_\A$.

\begin{defn}
Let $U, \, Y \in \RHinf$ and $\M:\RHinf \to \RHinf$ be of compatible dimension.  The optimization problem
\begin{equation}
\begin{array}{l}
\minimize{U\in \RHinf}  \Htwonorm{Y - \M(U)}^2+\rho\Htwonorm{U}^2+2\lambda\left\|U\right\|_{\A} \
\text{s.t. }  U \in \s
\end{array}
\label{eq:h2_opt}
\end{equation}
is called the $\Htwo$ \emph{RFD optimization problem} with parameters $(\rho,\lambda)$, distributed constraint $\s$, and atomic norm penalty $\left\|\cdot\right\|_{\A}$.
\end{defn}

There are two components of note in this definition.  The first is that $\rho$ need not be equal to $\rho_u$, the control cost parameter of the original non-penalized control problem \eqref{eq:the_opt}; the reasons why a different choice of $\rho$ may be desirable are explained in \S\ref{sec:recovery}.  
The second is the use of an atomic norm penalty function to induce structure.  Indeed, if one seeks a solution $U_*$ that can be composed as a linear combination of a small number of ``atoms'' $\A$, then a useful approach, as described in \ellcitations, \nuccitations, \cite{CRPW12}, to induce such structure in the solution of an inverse problem is to employ a convex penalty function that is given by the atomic norm induced by the atoms $\A$ \approxcitations. Examples of the types of structured solutions one may desire in linear inverse problems include sparse, group sparse and signed vectors,  and low-rank, permutation and orthogonal matrices (see \cite{CRPW12} for a more extensive list).

Specifically, if one assumes
that
\begin{equation}
U=\sum_{i=1}^{r}c_{i}A_{i},\ A_{i}\in\A,\ c_{i}\geq0
\label{eq:decomp}
\end{equation}
for a set of appropriately scaled and centered atoms $\A$,
and a small number $r$ relative to the ambient dimension, then solving
\begin{equation}
\begin{array}{rl}
\underset{U}{\mathrm{minimize}}  &  \cost{U;Y,\M} + 2\lambda\|U\|_\A
\end{array}
\end{equation}
with the atomic norm $\|\cdot\|_\A$ given by the gauge function\footnote{If no {such $t$ exists}, then $\|X\|_\A = \infty$.} 
\begin{equation}
\begin{array}{rcl}
||U||_{\A}: & = & \inf\{t\geq0\, \big{|}\, U\in t\mathrm{conv}(\A)\}\\
 & = & \inf\{\sum_{A\in\A}c_{A}\, \big{|}\, U=\sum_{A\in\A}c_{A}A,\, c_A \geq 0\}
\end{array}
\label{eqn:atomnorm}
\end{equation}
results in solutions that are both consistent with the data as measured in terms of the cost function $\cost{\cdot;Y,\Cl}$, and that are sparse in terms of their atomic descriptions, i.e., are a nonnegative combination of a small number of elements from $\A$.  \edit{Note that the unit ball of the atomic norm is given by the convex hull of its constituent atoms -- as such, atomic norms are convex functions.}



Our discussion in $\S$\ref{sec:RFDQI} on designing controller architecture by finding structured solutions to the model matching problem \eqref{eq:gen_mm} suggests natural atomic sets for constructing suitable penalty functions for RFD.  We make this point precise by showing that actuator, sensor, and/or communication delay design can all be performed through the use of a purposefully constructed atomic norm. 
 We introduce several novel penalty functions for controller architecture design, most notably for the simultaneous design of actuator, sensor and communication delays.   Further, all regularizers that have been considered for control architecture design in the literature (cf. \cite{LFJ13,JRM14,FLJ13,MC_CDC14,DJL14,M_CDC13_codesign, M_TCNS14}, among others) may be viewed as special instances of the atomic norms described below.

%
In what follows, the atomic sets that we define are of the form
\begin{equation}
\A = \displaystyle \bigcup_{\Aa \in \mathscr{M}} \Aa \cap k_{\Aa} \mathcal{B}_{\Htwo},
\label{eq:asets}
\end{equation}
for $\mathscr{M}$ an appropriate set of subspaces, $\{k_{\Aa}\}$ a set of normalization constants indexed by the subspaces $\Aa \in \m$, and $\mathcal{B}_{\Htwo}$ the $\Htwo$ unit norm ball; see the concrete examples below.  
Note that we normalize our atoms relative to the $\Htwo$ norm as this norm is isotropic; hence this normalization ensures that no atom is preferred over another within a given family of atoms $\A$.  We use $n_s$ and $n_a$ to denote the total number of sensors and actuators, respectively, available for the RFD task.

\edit{
\subsubsection{Actuator/Sensor Norm} \note{For the \emph{Actuator Norm}, we choose the atomic set to be transfer functions in $\RHinf^{n_a \times n_s}$ that have exactly one nonzero row with unit $\Htwo$ norm, i.e., suitably normalized Youla parameters that use only one actuator.  Specifically, the set of subspaces \eqref{eq:asets} in this context is
\begin{equation}
\m_{\text{act}} := \left\{ \Aa \subset \RHinf^{n_a \times n_s} \, \big{|} \, \text{$\Aa$ has one nonzero row} \right\},
\label{eq:act_sub}
\end{equation}
 leading to the atomic set
\begin{equation}
\A_{\text{act}} := \left\{ e_i V \, \big{|} \, V \in \RHinf^{1 \times n_s},\, \Htwonorm{V} = 1\right\}.
\label{eq:act_atoms}
\end{equation}
The resulting atomic norm is then given by
\begin{equation}
\norm{U}_{\text{act}} = \sum_{i=1}^{n_a} \Htwonorm{e_i^\top U}.
\label{eq:act_norm}
\end{equation}  In particular, each ``group'' corresponds to a row of the Youla parameter.  For the \emph{Sensor Norm}, we similarly choose transfer functions with exactly one nonzero column with unit $\Htwo$ norm, leading to the atomic norm }
\begin{equation}
\norm{U}_{\text{sns}} = \sum_{i=1}^{n_s} \Htwonorm{Ue_i}.
\label{eq:sns_norm}
\end{equation}
Both of these norms are akin to a Group Lasso penalty \cite{GroupLasso}.}


\subsubsection{Joint Actuator and Sensor Norm}  Conceptually, each atom corresponds to a controller that uses only a small subset of actuators and sensors. As each row of the Youla parameter $U$ corresponds to an actuator and each column to a sensor, the atomic transfer matrices have support defined by a submatrix of $U(z)$.  Specifically, we choose atoms with at most $k_a$ actuators and $k_s$ sensors:
\begin{equation}
\begin{array}{l}
\m_{\text{act+sns}} := \left\{ \Aa \subset \RHinf^{n_s\times n_a} \, \big{|} \, \supp{\Aa} \text{is a submatrix} \right. \\ \left. \text{with at most $k_a$ nonzero rows}  \text{ and $k_s$ nonzero columns} \right\}
\end{array}
\label{eq:act+sns_atoms}
\end{equation}
The scaling terms $k_{\Aa}$ in the definition of the atomic set \eqref{eq:asets} are given by $k_{\Aa}= \left(\text{card}(\Aa)+.1\right)^{-\frac{1}{2}}$, and are necessary as some of the atoms are nested within others -- the additional $.1$ can be any positive constant, and controls how much an atom of larger cardinality is preferred over several atoms of lower cardinality.
The resulting \emph{Actuator+Sensor Norm} is then constructed according to \eqref{eq:asets} and
 is akin to the latent Group Lasso \cite{OJV11}.

\subsubsection{Communication Link Norm}
\note{As described in \S\ref{sec:RFDQI}, the communication design task is to select which additional links to introduce into an existing base communication graph.  An atom in $\A_{\text{comm}}$ corresponds to such an additional link.  We provide an example of such an atomic set for a simple system, and refer the reader to \cite{M_TCNS14} for a more general construction.  
\begin{figure}[!h]
\centering
\includegraphics[width = .25\textwidth]{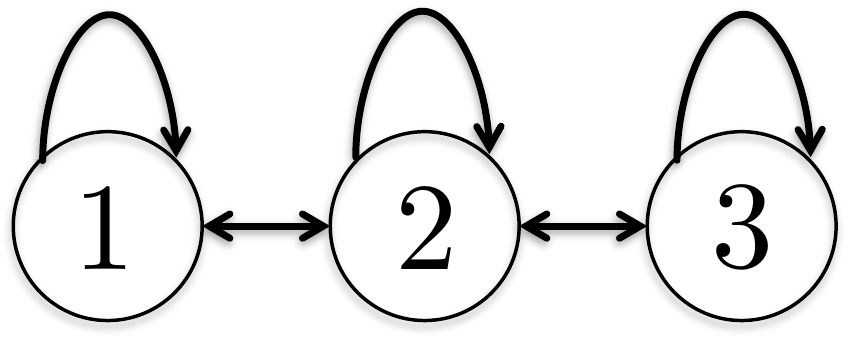}
\caption{Three player chain system}
\label{fig:4chain}
\end{figure}
In particular, consider a three player chain system, with physical topology illustrated in Figure \ref{fig:4chain}, such that $P_{22}$ lies in the subspace
\[
\s:= \frac{1}{z}\begin{bmatrix} \ast & 0 & 0  \\ 0 & \ast & 0 \\ 0 & 0 & \ast \end{bmatrix} \oplus 
\frac{1}{z^2}\begin{bmatrix} \ast & \ast & 0  \\  \ast & \ast & \ast \\ 0 & \ast & \ast \end{bmatrix} \oplus 
\frac{1}{z^3}\Rp,
\]
where $\ast$ is used to denote $\R$ to reduce notational clutter.
We consider an existing communication graph matching the physical topology illustrated in Figure \ref{fig:4chain} so that the induced distributed subspace constraint, as described in \S\ref{sec:RFDQI}, is given precisely by $\s$.  It can be checked that $\s$ is then QI with respect to $P_{22}$.  We consider choosing from two additional links to augment the communication graph: a directed link from node 1 to node 3, and a directed link from node 3 to node 1.  Then $\m_{\text{comm}} = \{\Aa_{13} , \Aa_{31} $\}, where its component subspaces are given by
\begin{equation*}
\begin{array}{rclrcl}
\Aa_{13} &=& \frac{1}{z^2}\begin{bmatrix} 0 & 0 & 0  \\  0 & 0 & 0 \\ \ast & 0 & 0 \end{bmatrix}, 
\Aa_{31} &=& \frac{1}{z^2}\begin{bmatrix} 0 & 0 & \ast  \\  0 & 0 & 0 \\ 0 & 0 & 0 \end{bmatrix}.
\end{array}
\end{equation*}
In particular, each subspace $\Aa_{ij}$ corresponds to the additional information available to the controller \emph{uniquely} due to the added link from sensor $j$ to actuator $i$.
The resulting \emph{Communication Link Norm} $\norm{\cdot}_{\text{comm}}$ is then constructed according to \eqref{eq:asets} with all normalization constants $k_{\Aa}=1$.
We note that this penalty is also akin to the latent Group Lasso \cite{OJV11}.
}


\subsubsection{Joint Actuator (and/or Sensor) and Communication Link Norm}

This penalty can be viewed as simultaneously inducing sparsity at the communication link level, while further inducing row sparsity as well.  The general strategy is to combine the actuator and communication link penalties in a convex manner.  We suggest two such approaches, one based on taking their weighted sums and the other based on taking their ``weighted maximum.''  In particular, we define the joint actuator \emph{plus} communication link penalty to be:
\begin{equation}
\norm{U}_{\text{act+comm}} = (1-\theta)\norm{U}_{\text{comm}} + \theta\norm{U}_{\text{act}},
\end{equation}
for some $\theta\in [0,1]$,
and the \emph{max} actuator/communication link penalty to be
\begin{equation}
\norm{U}_{\max\{\text{act},\text{comm}\}} = \max\left\{(1-\theta)\norm{U}_{\text{comm}}, \theta\norm{U}_{\text{act}}\right\},
\end{equation}
for some $\theta \in [0,1]$.
The analogous Sensor and Communication Link penalties, as well as Sensor+Actuator and Communication link penalties can be derived by replacing $\A_{\text{act}}$ with either $\A_{\text{sns}}$ or $\A_{\text{act+sns}}$.

\subsection{Further Connections with Structured Inference}
As already noted in \S\ref{sec:cost_fcns}, by choosing different values of $\rho$ for $\lambda=0$ we are able to recover control-theoretic analogs to Ordinary Least Squares and Ridge Regression \cite{Ridge}.  Noting that the actuator norm penalty \eqref{eq:act_norm} is akin to the Group Lasso \cite{GroupLasso}, we now discuss how control theoretic analogs of the Group Lasso and Group Elastic Net \cite{ElasticNet}  can be obtained by setting $\lambda >0$ in \eqref{eq:h2_opt} and using the Actuator Norm \eqref{eq:act_norm} penalty -- these connections are summarized in Table \ref{tab:dict}.    To simplify the discussion, we once again consider these connections in the context of the basic LQR problem introduced in Example \ref{ex:bLQR}, now augmented with the actuator norm penalty \eqref{eq:act_norm}.

\begin{table*}[t]
\scriptsize
\centering
\begin{tabularx}{\textwidth}{|S|SSs|SSs|}
\hline
 Parameters & Stats Problem & Stats Structure & Stats Tradeoff &  RFD Problem & {RFD~Structure} & RFD Tradeoff  \\ \hline
 $\lambda = 0$, $\rho=0$ &Ordinary Least Squares & None & N/A & Cheap Control  LQR & Deadbeat response & N/A   \\ \hline
 $\lambda = 0$, $\rho>0$ &Ridge Regression & Small Euclidean norm  & Bias, Variance & LQR & Small control action& State deviation, Control effort \\ \hline
 $\lambda > 0$, $\rho=0$ & Group LASSO & Group sparsity & Bias, Variance,  Model complexity & RFD LQR &  Sparse actuation & State deviation, Control effort, Actuation complexity \\ \hline
  $\lambda > 0$, $\rho>0$ & Group Elastic Net & Correlated group sparsity & Bias, Variance,  Model complexity & RFD LQR & Correlated  sparse actuation & State deviation, Control effort,  Actuation complexity \\ \hline
\end{tabularx}
\caption{A dictionary relating various SLIP problems in statistical inference and Actuator RFD problems.}
\label{tab:dict}
\end{table*}

%

\note{
In structured inference problems, the setting $\lambda>0,\, \rho >0$ corresponds to Group Elastic-Net regression.  If the groups are single elements, this becomes the traditional Elastic Net and Lasso.  The singleton group setting with $\lambda>0$, $\rho = 0$ corresponds to Lasso regression, and this inference method is employed when the underlying model is known to be sparse -- in particular, the Lasso penalty is used to select which elements $U^*_i$ of the model are non-zero \ellcitations.  Continuing with the singleton group setting, if both $\lambda>0$ and $\rho>0$, then the corresponding inferential approach is called the Elastic Net.  In addition to the sparsity-inducing properties of the Lasso, the Elastic Net also encourages automatic clustering of the elements \cite{ElasticNet} -- in particular, $\rho>0$ encourages the simultaneous selection of highly correlated elements (two elements $U^*_i$ and $U^*_j$ are said to be highly correlated if $\Cl(U^*_i) \approx \Cl(U^*_j)$). Thus $\rho$ can be seen as a parameter that can be adjusted to leverage a prior of \emph{correlation} in the underlying measurement operator $\M$.  These interpretations carry over to more general groups in a natural way.
}

In RFD, this setting corresponds to our motivating Example \ref{ex:bLQR} augmented with the actuator norm penalty, in which we design the controller's actuation architecture.  As each atom corresponds to an actuator, this RFD procedure then selects a small number of actuators.  Porting the clustering effect interpretation from the structured inference setting, we see that $\rho$ promotes the selection of actuators that have similar effects on the closed loop response.  In particular, this suggests that for systems in which no such similarities are expected, $\rho$ should be chosen to be small (or 0) during the RFD process, even if the original LQR problem had non-zero control cost.

\squeezeup\squeezeup
\edit{
\section{The RFD Procedure}
\label{sec:rfd_proc}
\subsection{The Two-Step Algorithm}
We now introduce the convex optimization based RFD procedure for the co-design of an optimal controller and the architecture on which it is implemented.  The remaining computational challenge is the possibly infinite dimensional nature of the RFD optimization problem \eqref{eq:rfd_opt}.  To address this issue, we propose a two step procedure: first, a finite dimensional approximation of optimization problem \eqref{eq:rfd_opt} is solved to identify a potential controller architecture and its defining subspace constraint $\D$.  Once this architecture has been identified, a traditional (and possibly infinite dimensional) optimal control problem \eqref{eq:trad_opt} with Youla prameter restricted to lie in $\D\cap\s$ is then solved -- in particular, in many interesting settings the resulting optimal controller restricted to the designed architecture can then be computed exactly leveraging results from the optimal controller synthesis literature \distcitations.

Formally, we begin by fixing an optimization horizon $t$ and a controller order $v$.  We suggest initially choosing $t$ and $v$ to be small (i.e., 2 or 3), and then gradually increasing these parameters until a suitable controller architecture/control law pair is found.  Our motivations for this approach are twofold: (i) first, selecting a small horizon $t$ and small controller order $v$ leads to a smaller optimization problem that is computationally easier to solve; and (ii) as we show in the next section, a smaller horizon $t$ and smaller controller order $v$ can actually aid in the identification of optimal controller architectures.  For a given performance metric $\cost{\cdot;Y,\Cl}$ and atomic norm penalty $\|\cdot\|_\Aa$, the two step RFD procedure consists of an architecture design step and an optimal control law design step:

\textbf{1) Architecture design:} Select the regularization weight $\lambda$ and solve the \emph{finite dimensional} RFD optimization problem
\begin{equation}
\begin{array}{rl}
\minimize{U \in \RHinf^{\leq v} \cap \s} & \cost{U;Y^{\leq t}, \Cl^{\leq t,v}} + 2\lambda \left\|U\right\|_\Aa
\end{array}
\label{eq:rfd_opt_trunc1}
\end{equation}
The actuators, sensors and communication links defining the designed architecture are specified by the non-zero atoms that constitute the solution $\hat{U}$ to optimization problem \eqref{eq:rfd_opt_trunc1}.  The architectural components employed in $\hat{U}$ in turn define a subspace $\D(\hat{U})$ which corresponds to all controllers (within the Youla parameterization) that have the same architecture as $\hat{U}$.\footnote{\edit{As described in \S\ref{sec:RFD}, the subspace $\D(\hat{U}) \cap \s$ is QI by construction, and hence this subspace also corresponds to all controllers with the same architecture as $\hat{K}= (I+\hat{U}P_{22})^{-1}\hat{U}$.}}

\textbf{2) Optimal control law design:} Solve the infinite dimensional optimal control problem with Youla parameter additionally constrained to lie in the designed subspace $\D(\hat{U})\cap \s$:
\begin{equation}
\begin{array}{rl}
\minimize{U \in \RHinf} & \cost{U;Y, \Cl} \
\text{s.t. } U \in \D(\hat{U}) \cap \s.
\end{array}
\label{eq:rfd_opt_full}
\end{equation}
If the resulting controller architecture and controller performance are acceptable, the RFD procedure terminates.  Otherwise, adjust $\lambda$ accordingly to vary the tradeoff between architectural complexity and closed loop performance.  If no suitable controller architecture/controller can be found, increase $t$ and $v$ and repeat the procedure. 
\begin{rem}[Removing Bias] \note{The method of solving a regularized optimization problem to identify the architectural structure of a controller and then solving a standard model matching problem restricted to the identified architecture is analogous to a procedure that is commonly employed in structured inference.  In structured inference problems, a regularized problem is solved first to identify a subspace corresponding to the structure of an underlying model $U_*$.  Subsequently, a non-regularized optimization problem with solution restricted to that identified subspace is solved to obtain an unbiased estimator of the underlying model $U_*$.}
\end{rem}

\squeezeup
\subsection{Simultaneous Actuator, Sensor and Communication RFD}
\label{sec:large}

In this subsection we demonstrate the full power and flexibility of the RFD framework in designing a distributed controller architecture, jointly incorporating actuator, sensor and communication link design.  In particular we consider a plant with eleven subsystems with topology as illustrated in Figure \ref{fig:RFD_topo}.  The solid lines correspond to the physical interconnection between subsystems. Choosing $C_2 = B_2 = I$, the adjacency matrix of this graph then defines the support of the $A$ matrix in the state space realization 
of the generalized plant \eqref{eq:genplant}, as well as the required communication links between nodes such that the distributed constraint is QI under $P_{22}$ \cite{M_TCNS14,RCL10}. 


\begin{figure}[h!]
\centering
\includegraphics[width = .35\textwidth]{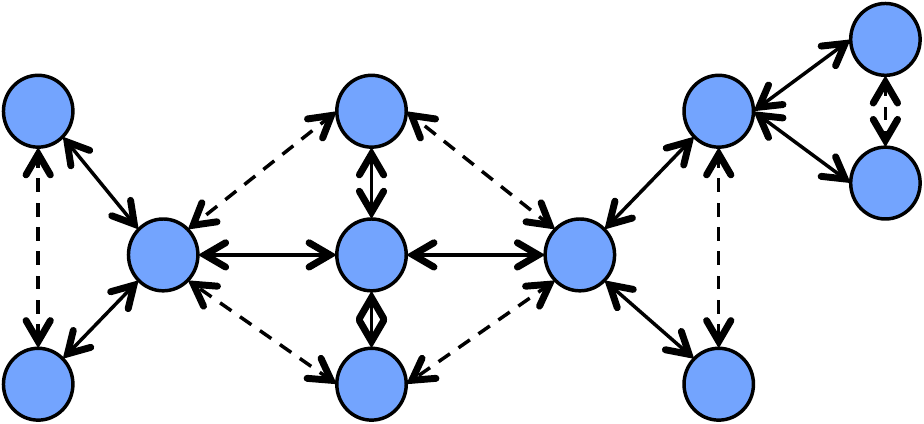}
\caption{Topology of system considered for RFD example.  Solid lines indicate both physical interconnections and existing communication links between controllers.  Dashed lines correspond to possible additional edges to be added.  
}
\label{fig:RFD_topo}
\end{figure}
\vspace{-1mm}
The non-zero entries of $A$ were generated randomly and normalized such that $|\lambda_{\max}\left(A\right)|=.999$. The remaining state space parameters of the generalized plant \eqref{eq:genplant} satisfy $D_{12}\begin{bmatrix} C_1^\top & D_{12}^\top \end{bmatrix} =  \begin{bmatrix} 0 & 25 I \end{bmatrix}$, and $D_{21}^\top\begin{bmatrix} B_1 & D_{21} \end{bmatrix} =  \begin{bmatrix} 0 & .01 I \end{bmatrix}$, with $C_1C_1^\top = 100 I$ and $B_1^\top B_1 = I$.

For the RFD task, we choose an $\Htwo$ norm performance metric; we allow each node to be equipped with an actuator and/or a sensor (for a total of 11 possible actuators and sensors), and we allow the communication graph to be augmented with any subset of the interconnections denoted by the dashed lines in Figure \ref{fig:RFD_topo}, in addition to the already present links given by the solid lines.  This leads to 536,870,911 different possible controller architectures.

We solved the RFD optimization \eqref{eq:output} with atomic norm $\norm{\cdot}_{\text{act+sns+comm}}$ as defined in \S\ref{sec:atoms}, with weighting parameter $\theta = .75$ and with $k_s = k_a = 1$.  \edit{We performed the RFD procedure for two different horizon/order pairs:  $t=4$ and $v=2$, as well as $t=6$ and $v=3$; for these latter horizon/order values acceptable tradeoffs between architecture complexity and closed loop performance were identified, and hence the RFD procedure terminated.  For each horizon$/$order pair $(t,v)$, we vary $\lambda$, and for each resulting optimal solution $\hat{U}$, we identified the designed architecture and corresponding subspace $\D(\hat{U})$.  We then used the method from \cite{LD14} to exactly solve the resulting non-regularized distributed $\Htwo$ model matching problem with subspace constraint $\D(\hat{U})$.}  Note that the Youla parameter solving this non-regularized problem is not restricted to have a finite impulse response.  In particular, we can compute the optimal Youla parameter and the corresponding optimal controller restricted to the architecture underlying $\hat{U}$ in a computationally tractable fashion because we guarantee that the subspace corresponding to the designed architecture is QI, as per the discussion in Section \ref{sec:RFD_algos}.

For horizon $t=6$ and order $v=3$, the resulting architectural complexity is plotted against the closed loop norm of the system in Figure \ref{fig:rfd_perf}.  As $\lambda$ is increased, the architectural complexity (i.e. the number of actuators, sensors and communication links) decreases, but at the expense of deviations from the performance achieved by the controller that uses all of the available architectural resources.  We also show the resulting architecture for $\lambda = 500$ in Figure \ref{fig:designed}: as can be seen, a non-obvious combination of eight actuators, eight sensors and five additional communication links are chosen, resulting in only a 0.71\% degradation in performance over the distributed controller using all eleven actuators, eleven sensors and seven additional communication links.

\begin{figure}[h!]
\centering
\includegraphics[width = .5\textwidth]{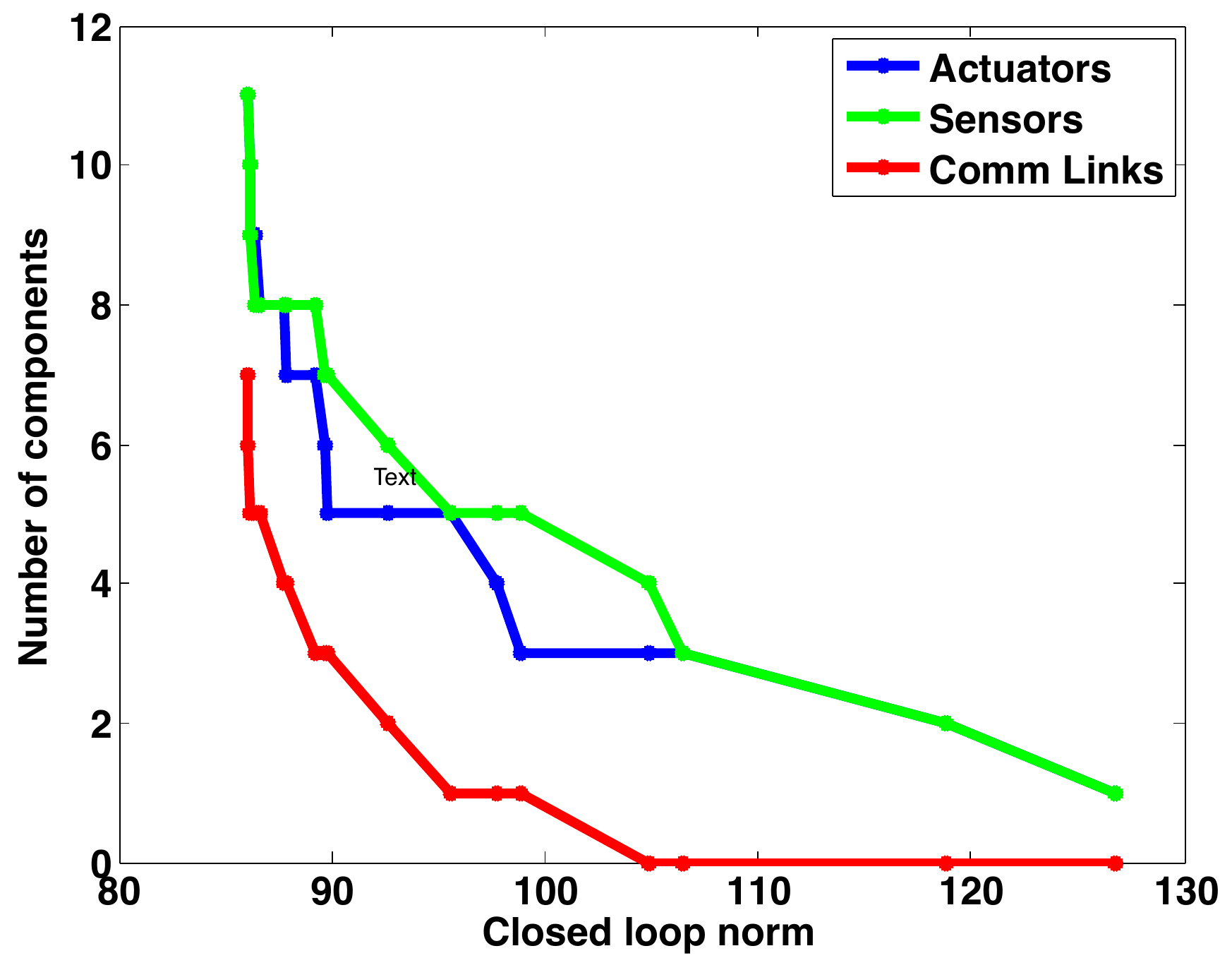}
\caption{A small degradation in closed loop performance allows for a significant decrease in architectural complexity.}
\label{fig:rfd_perf}
\end{figure}


\begin{figure}[h!]
\centering
\includegraphics[width = .35\textwidth]{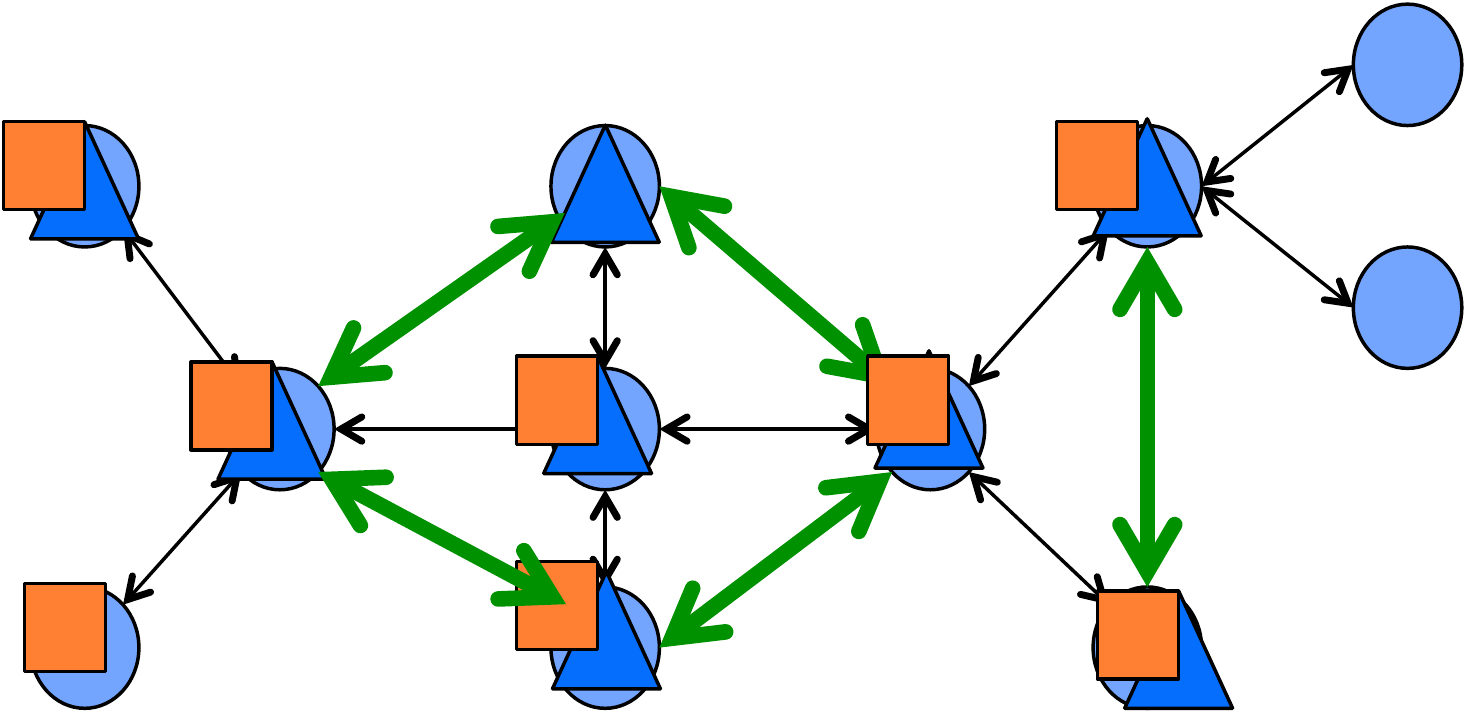}
\caption{Resulting architecture for $\lambda=500$: despite only using eight actuators (orange squares), eight sensors (blue triangles) and five additional communication links (green arrows), the performance only degraded by 0.71\% relative to  the distributed controller using all eleven actuators, eleven sensors and seven communication links.}
\label{fig:designed}
\end{figure}

As this example shows, the RFD procedure is effective at identifying simple controller architectures that approximate the performance of a controller that maximally utilizes the available architectural resources. In the next section, we offer some theoretical justification for the success of our procedure by suitably interpreting the RFD optimization problem \eqref{eq:rfd_opt} in the context of approximation theory and by making connections to analogous problems in structured inference.
}

\section{Recovery of Optimal Actuation Structure}
\label{sec:recovery}
This section is dedicated to the analysis of the $\Htwo$ RFD optimization problem \eqref{eq:h2_opt} with no distributed constraint, actuator norm penalty \eqref{eq:act_norm}, and $Y$ and $\M$ as given in Example \ref{ex:actuators} -- a nearly identical argument applies to a sensor norm regularized problem.  We discuss how to extend the analysis to output feedback and distributed problems at the end of the section. 

\edit{
\note{Viewing the model matching problem \eqref{eq:gen_mm} as a linear inverse problem makes it clear that designing a structured controller is akin to obtaining a structured solution to a linear inverse problem.  The problem of obtaining structured solutions to linear inverse problems arises prominently in many contexts, most notably in statistical estimation.  In that setting, one posits}
%
 a linear measurement model\footnote{We purposefully use non-standard notation to facilitate comparisons between RFD and SLIP.}
\begin{equation}
Y = \Cl(U_*) + W,
\label{eq:truth}
\end{equation}
 where $Y$ is the vector of observations, $\Cl$ is the measurement map, $U_*$ parametrizes an underlying model and $W$ is the measurement error.  
\note{The linear model \eqref{eq:truth} also has an appealing interpretation from a control-theoretic perspective.
 In particular, letting $Y\in\RHinf$ be the state component of the open loop response of a LTI system, $U_*\in \RHinf$ be a suitably defined Youla parameter, and $\Cl:\RHinf\to\RHinf$ be the map from Youla parameter to the state component of the closed loop response,  it is then immediate that
\begin{equation}
 W := Y -\Cl(U_*)
\label{eq:opt_ctrl}
\end{equation}
represents the state component of the closed loop response achieved by the controller $U_*$.  Table \ref{tab:pm} summarizes the correspondence between these two perspectives.}

\begin{table}[h]
\centering
\begin{tabular}{ccc}
\hline
 Parameter &  Structured Controller Design & Structured Inference  \\ \hline
 $Y$ &Open loop system  & Observations   \\ \hline
 $\Cl$ &Map to closed loop  & Measurement map   \\ \hline
 $U_*$ &Desired controller & Underlying model  \\ \hline
 $W$ & Closed loop response & Measurement noise \\ \hline
\end{tabular}
\caption{Interpretation of parameters in Structured Controller Design and Structured Inference. }
\label{tab:pm}
\end{table}

This conceptual connection suggests a novel interpretation of the role of the closed loop response $W$ achieved by a controller $U_*$. In an inferential context, since $W$ corresponds to measurement noise, a smaller $W$ makes the task of identifying the structure of the underlying model $U_*$ much easier, as the measurements are more accurate.  In a similar spirit, we demonstrate that structured controller design is easier (via the solution to an RFD optimization problem \eqref{eq:rfd_opt}) if the corresponding \note{state component} of the closed loop response is small. Thus
the state component of the closed loop response of the system plays the role of noise when trying to identify the structure of a suitably defined controller $U_*$.  In the sequel, we describe an appropriate notion of smallness for the \note{state component} of the closed loop response in the context of designing structured controllers.

}

\edit{
The remainder of the discussion in this section builds on prominent results from the structured inference literature \cite{CRT06,Donoho06,fazelThesis,RFP10}.  The flavor of these results is somewhat non-standard in the controls literature, and we therefore pause briefly to frame the setup in this section appropriately and to discuss how the results of this section should be interpreted.  The main result of this section proceeds by assuming that there exists an architecturally simple controller $U_*$  (i.e., one with a small number of actuators) that achieves a good closed loop response, (i.e., that achieves a small $W$ as defined in \eqref{eq:opt_ctrl}).  Under suitable conditions, Theorems \ref{thm:consistency} and \ref{thm:main} state that the architectural structure of $U_*$ can be recovered via tractable convex optimization using the RFD procedure.  These conditions are phrased in terms of quantities associated to $U_*$ which are typically unknown in advance -- however, these conditions are not meant to be checked prior to solving a RFD optimization problem.  Although the results are stated in terms of a nominal controller $U_*$, they should be interpreted as describing the properties satisfied by controller architectures identified via the RFD procedure of \S\ref{sec:rfd_proc}.  In particular, the RFD procedure requires solving RFD optimization problems across a range of controller orders $v$, optimization horizons $t$ and regularization weights $\lambda$: this process leads to a set of controller architectures being identified.  Our results allow a practitioner to be confident that all controller architectures satisfying the conditions of our theorems -- i.e., those that have a small number of actuators and that achieve a small closed loop state response -- are included in this set of identified controller architectures.  In this way, the RFD procedure provably identifies good controller architectures, should they exist.} 

We study finite dimensional variants of the $\Htwo$ RFD optimization problem \eqref{eq:h2_opt} with the actuator norm penalty \eqref{eq:act_norm}, and show that such finite dimensional approximations are sufficient to identify the {structure} of a desired controller $U_\ast$.  
\note{In particular, we truncate the optimization problem \eqref{eq:h2_opt} to a finite horizon $t$ by restricting $Y-\Cl(U)\in\RHinf$ to the first $t$ elements of its impulse response, and to a finite controller order $v$ by restricting $U$ to lie in $\RHinf^{\leq v}$.  The resulting optimization problem is thus finite dimensional, and \edit{corresponds to the first step of the RFD procedure defined in the previous section.}
}
At this point, it is convenient to introduce the temporally truncated version of \eqref{eq:opt_ctrl} for a fixed optimization horizon $t$ and controller order $v$:
\begin{equation}
\begin{array}{rcl}
W^{\leq t} &=& Y^{\leq t} - \M^{\leq t,t}(U^{\leq t}_*) \\
&=& Y^{\leq t}  - \M^{\leq t,v}(U^{\leq v}_*) - \drv
\end{array}
\label{eq:Ttruth}
\end{equation}
with \begin{equation}\drv:= \M^{\leq t,t}\left(U_*^{\leq t} - U_*^{\leq v}\right)\label{eq:tail}\end{equation}
corresponding to the effect of the ``tail'' of $U_*$ on the state component $W^{\leq t}$ of the truncated closed loop response.  

The flexibility in the choice of the optimization horizon $t$ and controller order $v$ will be the focus of much of our discussion.  In particular, it is of interest to find the smallest $t$ and $v$ for which we can guarantee that the RFD procedure recovers the structure underlying $U_*$ -- the smaller the horizon and controller order, the smaller the size of the optimization problem that needs to be solved.  Perhaps counter-intuitively, we show that larger $t$ and $v$ do not necessarily help in recovering the structure of an underlying parameter $U_*$.  We make this statement precise in what follows, but again drawing on intuition from the structured inference literature, we note that increasing $v$ in RFD is analogous to increasing the allowed model complexity when solving an inference problem: if the model class is too rich, we risk over-fitting and thus obfuscating the structure of the underlying model $U_*$.



Our goal is to prove that the solution $\tilde{U}$ to the finite dimensional $\Htwo$ RFD optimization problem 
\begin{equation}
\begin{array}{l}
\tilde{U} = \argmin{U \in {\RHinf^{\leq v}}}  \Htwonorm{{Y^{\leq t}+\Mtv(U)}}^2   
+\rho\Htwonorm{U}^2 +2\lambda\Gnorm{U} 
\end{array}
\label{eq:Topt}
\end{equation}
has the same architectural structure as $U_*$ for appropriately chosen $t$ and $v$.  To show this, we study the solution $\hat{U}$ to the following \emph{architect} optimization problem:
\begin{equation}
\begin{array}{l}
\hat{U} = \argmin{U \in \RHinf^{\leq v}}  \Htwonorm{{Y^{\leq t}+\Mtv(U)}}^2 +\rho\Htwonorm{U}^2 +2\lambda\Gnorm{U} \\
\indent \indent \indent \indent \, \text{s.t. }  U \in \Gs
\end{array}
\label{eq:Toracle}
\end{equation}
where $\Gs$ is a subspace of Youla parameters $U$ with the property that a row of $U_*$ being zero implies that the corresponding row of $U$ is zero.  In words, $\Gs$ may be viewed as the set of Youla parameters corresponding to actuation schemes matching the actuation scheme of $U_*$. We also define $\m_* \subset \m_{\text{act}}$, with $\m_{\text{act}}$ defined as in \eqref{eq:act_sub}, to be
\begin{equation}
\m_* :=  \left\{\Aa \in \m_{\text{act}} \, {|} \,  (U_*)_{\Aa} \neq 0 \right\}.
\end{equation}
\note{In words, the elements of $\m_*$ correspond to actuation schemes that use a single actuator, where these actuators are defined by the nonzero rows of the desired controller $U_*$.}

We show under suitable conditions on $t$, $v$, $U_*$ and $\Mtv$ that $\hat{U}=\tilde{U}$; that is to say that the architect solution $\hat{U}$ is also the unique optimal solution to the  RFD optimization problem \eqref{eq:Topt} without the additional constraint $U \in \Aa_*$.  \note{As a result, since $\hat{U}$ is constrained to lie in $\Aa_*$, the solution to the RFD optimization problem $\tilde{U}$ also lies in $\Aa_*$ and hence has the same architectural structure as $U_*$.} We emphasize that at no stage during the RFD procedure described in \S\ref{sec:rfd_proc} do we require knowledge of $\Aa_*$ and $\m_*$ -- the investigation of the architect problem \eqref{eq:Toracle} is only a theoretical tool used to prove structural recovery results.


%
%

\subsection{Identifiability Conditions in Control}

We begin by introducing two \emph{restricted gains} in terms of the subspace $\Gs$ and its orthogonal complement $\Gs^\perp$.    In order to do so, we introduce the dual norm to $\norm{\cdot}_{\text{act}}$, which is given by
\begin{equation}
\Gdnorm{U} = \max_{\Aa \in \m_\text{act}}\Htwonorm{U_\Aa}.
\end{equation}

These restricted gains are then

\begin{equation}
\begin{array}{rl}
\alpha^{\leq t,v} := \displaystyle\min_{\Delta}& \Gdnorm{\left(\adj{\Mtva} \Mtva + \rho I\right)(\Delta)} \\
\text{s.t.} & {\Gdnorm{\Delta}=1}, \, \Delta \in \Gs\cap {\RHinf^{\leq v}}
\end{array}
\label{eq:alpha}
\end{equation}
\begin{equation}
\begin{array}{rl}
\beta^{\leq t,v} := \max &  \Gdnorm{\adj{\Mtvap} \Mtva(\Delta)} \\
\text{s.t.} & {\Gdnorm{\Delta}\leq1}, \, \Delta \in \Gs\cap{\RHinf^{\leq v}}.
\end{array}
\label{eq:beta}
\end{equation}

The minimum gain $\atv$ is a quantitative measure of the injectivity of the operator $\Mtv\times \sqrt{\rho} I$ restricted to the subspace $\Gs$.  Intuitively, it characterizes the distinctions among the effects of the different actuators within $\Gs$.  The maximum gain $\btv$, on the other hand, is a measure of how different the effects of actuators in $\Gs$ are from those of actuators in $\Gs^\perp$.

We can already see some immediate implications of different choices of the horizon $t$ and and controller order $v$ on these quantities.  In particular, $\atv$ is non-increasing in the controller order $v$.  This minimum gain's dependence on the horizon $t$ is more subtle.  \edit{
Define the \emph{mixing time of $\m_*$} to be 
\begin{equation}
\mix := \max\left\{ t \in \mathbb{Z}_+ \, \big{|} \, \adj{\M^{\leq t}_{\Aa}}\M^{\leq t}_\Bb = 0, \ \forall \Aa \neq \Bb \in \m_*\right\}.
\label{eq:tau}
\end{equation}}
 If no $t$ exists such that the condition within the $\max\left\{\cdot\right\}$ is satisfied, we set $\mix = 0$.
 The mixing time $\mix$ measures how long it takes for the effects of the distinct actuators used by $U_*$ to overlap, or mix, in the closed loop  response.  Consequently the minimum gain $\atv$ is non-decreasing in $t$ so long as $t\leq \mix$, i.e., so long as $t$ is sufficiently small that the effects of the different actuators used by $U_*$ do not overlap.  
We then have the following lemma:
\begin{lem}
Let $\mix$ be as defined in \eqref{eq:tau}.  Then 
\begin{equation}
\atv = \rho + \min_{\Aa\in\m_*}\sm{\adj{\Mtv_\Aa}\Mtv_\Aa}
\label{eq:alpha_lem}
\end{equation}
for all $1\leq t \leq \mix$.  In particular, $\atv$ is non-decreasing in $t$ for all $1\leq t \leq \mix$.
\label{lem:alpha_smallT}
\end{lem}
\begin{IEEEproof}
It is easily verified that for $\Aa\neq \Bb \in \m_*$ and $t \leq \mix$, we have $\adj{\Mtv_\Aa}\Mtv_\Bb = 0$, from which \eqref{eq:alpha_lem} follows.  To see that $\atv$ as given in \eqref{eq:alpha_lem} is non-decreasing in $t$, it suffices to note that
\begin{equation}
\adj{\M_\Aa^{\leq t+1,v}}\M_\Aa^{\leq t,v} = \adj{\Mtv_\Aa}\Mtv_\Aa,
\end{equation}
leading to the conclusion that
\begin{equation}
\begin{array}{l}
\adj{\M_\Aa^{\leq t+1,v}}\M_\Aa^{\leq t+1,v} = \adj{\Mtv_\Aa}\Mtv_\Aa +  \adj{\M_\Aa^{\leq t+1,v}-\Mtv_\Aa}\left({\M_\Aa^{\leq t+1,v}}-\Mtv_\Aa\right).
\end{array}
\end{equation}
The result follows by noting that the final term in this expression is positive semidefinite.
\end{IEEEproof}

In particular, this result suggests that actuation schemes with more evenly distributed actuators (i.e., those with larger mixing times $\mix$ \eqref{eq:tau}) are easier to identify.

The maximum gain $\btv$, however, is clearly seen to be non-decreasing both in the controller order $v$ and the horizon $t$.  This is consistent with our interpretation of $\btv$ as a measure of similarity between actuators: as either $v$ or $t$ increase, there is more time for the mixing of the actuators' control actions via the propagation of dynamics in the system, increasing their worst-case ``similarity.''
We now assume that the following \emph{identifiability condition} is satisfied.
\begin{assumption}[Identifiability]
There exist $1\leq v \leq t < \infty$  such that
\begin{equation}
\frac{\btv}{\atv} =: \nu\in\left[0,1\right).
\label{eq:nu}
\end{equation}
\label{as:nu}
\end{assumption}

In light of the previous discussion, it is immediate that a larger controller order $v$ decreases the likelihood of the identifiability condition being satisfied, and should therefore be taken as small as possible.  The effect of increasing the horizon $t$ is less clear, but we see that it may help if the minimum gain $\atv$ increases sufficiently fast with $t$ relative to the increase in the gain $\btv$ with respect to $t$ -- further there is no need to increase $t$ beyond the mixing time $\mix$.

In the inference literature, the analog of these identifiability assumptions are given by conditions known as the \emph{restricted eigenvalue condition} \cite{RWY10} and the \emph{restricted isometry property} \cite{CT05}.  In the sequel, we give an example of {deterministic} and {structured} state space matrices that satisfy these identifiability conditions.  Specifically, we focus on systems \eqref{eq:genplant} that have block diagonal $B_2$ and $C_1$ matrices (i.e., decoupled actuators and state costs), and block banded state matrices $A$ (i.e., locally coupled dynamics).

\edit{
\begin{rem}
Notice that if $\Mtv = I$, then $\atv \geq 1$, $\btv = 0$ and $\nu = 0$.  These conditions are satisfied if $B=C=I$ and $v=t=1$ in Example \ref{ex:bLQR} (Basic LQR), or if $C_1 = B_2 = I$ and $ v =1$, $t = 2$ in Example \ref{ex:actuators} ($\Htwo$ State Feedback).  Thus sufficiently small values of $v$ and $t$ ensure that condition \eqref{eq:nu} holds.  However, the resulting optimization problem only incorporates low order effects of the dynamics (as encoded in $\Mtv$) in the RFD optimization problem, suggesting that $t$ and $v$ should be also be chosen large enough to sufficiently capture the dynamics of the system. This observation and Lemma \ref{lem:alpha_smallT} motivate our suggestion in Section \ref{sec:rfd_proc} to begin with small horizon $t$ and controller order $v$ and to then gradually increase these values until a suitable controller architecture/control law pair is found.
\end{rem}
}

\subsection{Sufficient Conditions for Recovery}




The following theorem provides sufficient conditions for (i) the architect solution $\hat{U}$ to be the unique optimal solution to the finite dimensional RFD optimization \eqref{eq:Topt}, and (ii) an actuator of the desired controller, identified by a subspace $\Aa \in \m_*$,  to be identified by the RFD procedure.


\begin{thm}[Structural Recovery]
Fix a horizon $1 \leq t < \infty$, and a controller order $1 \leq v \leq t$ such that Assumption \ref{as:nu} 
holds.
If
\begin{equation}
\begin{array}{l}
\lambda > \frac{\nu}{1-\nu}\left( \Gdnorm{\adj{\Mtv_{\Gs}} (\err)} + \rho \Gdnorm{{U_*^{\leq v}}}\right) +\frac{1}{1-\nu}{\Gdnorm{\adj{\Mtv_{\Gs^\perp}} (\err)}}
\end{array}
\label{eq:sufficient}
\end{equation}
we have that $\hat{U}$ as defined in \eqref{eq:Toracle} is the unique optimal solution to \eqref{eq:Topt}, and that the row support of ${\hat{U}}$ is contained within the row support of $U_*$.
Further if $\Aa \in \m_*$ and
\begin{equation}
\begin{array}{l}
\Htwonorm{{(U^{\leq v}_*)_\Aa}} > \frac{1}{\alpha}\left(\lambda + \Gdnorm{\adj{\Mtv_{\Gs}} (\err)} + \rho \Gdnorm{{U^{\leq v}_*}}\right),
\end{array}
\label{eq:structure}
\end{equation}
 then $\hat{U}_\Aa\neq 0$.
\label{thm:consistency}
\end{thm}

The condition \eqref{eq:sufficient}  states that, under suitable identifiability assumptions, the regularization weight $\lambda$ needs to be sufficiently large to guarantee that the architect solution $\hat{U}$ is also the solution to \eqref{eq:Toracle}.  However, this can always be made to hold by choosing $\lambda$ sufficiently large so that $\hat{U} = \tilde{U} = 0$.  The second condition \eqref{eq:structure}  provides an upper bound on the values of $\lambda$ for which a specific actuator (i.e., a specific component $(U^{\leq v}_*)_\Aa, \, \Aa \in \m_*$) is identified by the architect solution $\hat{U}$.  
%
%
The following corollary then guarantees the recovery of $\m_*$.
\begin{coro}
Let $\mu := \min_{\Aa \in \m_*} \|(U^{\leq v}_*)_\Aa\|_{\Htwo}$, and suppose that the open interval
\begin{equation}
\begin{array}{r}
\Lambda := \left(\frac{\nu}{1-\nu}\left[ \Gdnorm{\adj{\Mtv_{\Gs}} (\err)} + \rho \Gdnorm{{U_*^{\leq v}}}\right] +\frac{\Gdnorm{\adj{\Mtv_{\Gs^\perp}} (\err)}}{1-\nu}, \right.\\  
\indent\left.\atv\mu - \Gdnorm{\adj{\Mtv_{\Gs}} (\err)} - \rho \Gdnorm{{U^{\leq v}_*}} \right)
\end{array}
\label{eq:Lrange}
\end{equation}
is non-empty.  Then the solution $\tilde{U}$ to the RFD optimization \eqref{eq:Topt}, with any regularization weight $\lambda$ chosen within $\Lambda$, has row support \emph{equal} to that of $U_*$. 
\label{coro:success}
\end{coro}

\edit{Note that it is useful to have a given architecture be identifiable for a range of regularization weights $\lambda$, as prior information about the values needed to specify \eqref{eq:Lrange} are typically not available.  We exploited this fact when we defined the RFD procedure in Section \ref{sec:rfd_proc} by suggesting that $\lambda$ be varied until a suitable architecture/control law pair is identified.}  
This corollary \edit{also} makes explicit that larger values of $\rho$ shrinks the range of $\lambda$ for which the RFD procedure is successful.  It also shows that larger $\tail$ tail terms \eqref{eq:tail} are deleterious to the performance of the RFD procedure as well -- therefore although we previously stated that the controller order $v$ should be chosen as small as possible, it should not be so small that the tail term $\tail$ is too large.

\begin{rem}[Extension to Output Feedback]
\label{rem:output_feedback}
A similar argument applies to the output feedback problem, but at the expense of more complicated formulas.  In particular the $\Htwo$ RFD optimization takes the form
\begin{equation}
\begin{array}{rl}
\minimize{{U}\in\RHinf} & \Htwonorm{Y - \M({U})}^2 + \Htwonorm{\frak{F}(U)}^2+2\lambda\norm{{U}}_{\text{act}},
\end{array}
\label{eq:output}
\end{equation}
for $Y=P_{11}$, $\M(U) = P_{12}UP_{21}$, and $\frak{F}(U) = \begin{bmatrix} P_{12}UD_{21} & D_{12}UP_{21} & D_{12}UD_{21} \end{bmatrix}$.  \edit{Notice that if $D_{12}$ and $D_{21}$ are set to 0 in the RFD optimization problem \eqref{eq:output}, we recover an optimization problem of exactly the same form as \eqref{eq:h2_opt} with $\rho = 0$, in which case the analyses of this section and the next section are applicable.}
\end{rem}
\begin{rem}[Extension to Distributed Constraints]
For the purposes of analysis, the additional constraint that $U \in \s$ can be incorporated by considering the restriction of $\M$ to $\s$, resulting in a centralized problem (cf. \cite{RL06} for an example of how this can be done). 
\end{rem}

\subsection{A RFD Signal to Noise Ratio}
\label{sec:parameters}
Theorem \ref{thm:consistency} as stated does not yet provide an immediate interpretation of the effect of the choices of the horizon $t$ and the controller order $v$ on the performance of the RFD procedure.  
In order to better understand the effects of the horizon $t$ and controller order $v$ on the success of the RFD procedure, we describe more interpretable bounds on $\atv$ and $\btv$.

\begin{lem}
Fix $1\leq t<\infty$ and $1 \leq v \leq t$.  The parameters $\atv$ and $\btv$, as defined in equations \eqref{eq:alpha} and \eqref{eq:beta}, can be bounded from below and above, respectively, as follows:
\begin{equation}
\atv \geq \rho + \gtv,
\label{eq:alpha_bound}
\end{equation}
where
\begin{equation}
\begin{array}{l}
 \gtv := \displaystyle \min_{\B \subseteq \m_*, |\B|\geq 1}\max_{\Aa\in \B}  \left[\sigma_{\min}\left(\adj{\M_\Aa^{{\leq} t,v}} \Mtv_\Aa\right) 
 - \sum_{\Bb\neq \Aa \in \B} \sigma_{\max}\left(\adj{\Mtv_\Aa} \Mtv_\Bb \right)\right]
 \end{array}
\label{eq:gamma}
\end{equation}
and
\begin{equation}
\btv \leq \max_{\Aa \in (\m_{\text{act}}\backslash \m_*)} \sum_{\Bb\in\m_*}\sigma_{\max}\left(\adj{\Mtv_\Aa}\Mtv_\Bb \right).  
\label{eq:beta_bound}
\end{equation}
Consequently, we can upper bound the ratio \eqref{eq:nu} as
\begin{equation}
\nu \leq \frac{\btv}{\rho + \gtv}.
\label{eq:nu_bound}
\end{equation}
\label{lem:param_bounds}
\end{lem}

In particular, it is a straightforward consequence of Lemma \ref{lem:alpha_smallT} that  for all $t \leq \mix$ (where the mixing time $\mix$ is as in \eqref{eq:tau}), the intermediate quantity $\gtv$, as introduced in Lemma \ref{lem:param_bounds}, is given by
\[\gtv= \min_{\Aa\in \m_*}\sm{\adj{\Mtv_\Aa}\Mtv_\Aa},\]
and is non-decreasing in $t$.  Further it is easily verified that the bounds \eqref{eq:alpha_bound}, \eqref{eq:beta_bound} and \eqref{eq:nu_bound} can be taken with equality if $v=1$, as each $\M^{\leq t,1}_\Aa$ is isomorphic to a column vector.

\note{The bounds computed in Lemma \ref{lem:param_bounds} can be combined with the sufficient conditions of Theorem \ref{thm:consistency} to describe sufficient conditions for the successful recovery of the architecture of $U_*$ in terms of a signal to noise like quantity -- to that end, we introduce the following definitions.}

\begin{defn}[RFD Noise]
We define the \emph{RFD Noise} level $\Noise$ for an $\Htwo$ RFD optimization problem \eqref{eq:Topt} to be
\begin{equation}
\begin{array}{l}
\Noise := \Gdnorm{\adj{\Mtv_{\Gs}} (\err)} + \Gdnorm{\adj{\Mtv_{\Gs^\perp}} (\err)}.
\end{array}
\label{eq:noise}
\end{equation}
\label{def:noise}
\end{defn}

\note{The control theoretic interpretation \eqref{eq:opt_ctrl} of the linear model \eqref{eq:truth} used in inference problems motivates our terminology -- recall in particular that in \eqref{eq:opt_ctrl} the state component of the closed loop response $W$ is interpreted as measurement noise in the context of identifying a structured controller.  Likewise, $T^{\leq t,v}$ can be viewed as additional noise introduced into the architecture identification procedure by the temporal truncation procedure described in \eqref{eq:Ttruth}.  We proceed to define a control theoretic analog to the signal in the context of RFD optimization problems.}

\begin{defn}[RFD SNR]
In the context of architecture recovery via RFD, the magnitude of each atom, $\Htwonorm{(U^{\leq v}_*)_\Aa}$ plays the role of a signal, and the RFD Noise level $\Noise$ that of noise, leading to the definition of the SNR of a component $(U^{\leq v}_*)_\Aa, \, \Aa \in \m_*$ as
\begin{equation}
\SNR{{(U^{\leq v}_*)_\Aa}} := \frac{\Htwonorm{{(U^{\leq v}_*)_\Aa}}}{\Noise}.
\label{eq:groupSNR}
\end{equation}
\label{def:snr}
\end{defn}

These definitions allow us to state simple conditions in terms of the SNR \eqref{eq:groupSNR} for the successful recovery of an actuation architecture via the solution to the $\Htwo$ RFD optimization problem \eqref{eq:Topt}.



\edit{
\begin{thm}
Let $\rho=0$, $\lambda = \lambda' + \kappa$, where $\lambda'$ is given by the right hand side of \eqref{eq:sufficient}, and $\kappa>0$ is an arbitrarily small constant, and assume that $\nicefrac{\btv}{\gtv} < 1$.  If
\begin{equation}
\SNR{{(U^{\leq v}_*)_\Aa}} > \frac{1}{\gtv - \btv}
\label{eq:main}
\end{equation}
for all $\Aa \in \m_*$, then for sufficiently small $\kappa$, the solution $\tilde{U}$ to the $\Htwo$ RFD optimization problem \eqref{eq:Topt} has the same row support as $U^{\leq v}_*$.
\label{thm:main}
\end{thm}
\begin{IEEEproof}
Follows from rearranging terms in \eqref{eq:structure}, Definition \ref{def:snr} and letting $\kappa$ tend to 0 from above.
\end{IEEEproof}}

\note{Setting $\rho=0$ increases the range $\Lambda$, as defined in \eqref{eq:Lrange}, for which the RFD optimization problem is successful in recovering the structure of $U_*$, and the assumption that $\nicefrac{\btv}{\gtv} < 1$ ensures that Assumption \ref{as:nu} holds.  Thus Theorem \ref{thm:main} can be viewed as a slightly stronger, but more interpretable, set of sufficient conditions for the success of the RFD procedure.}

Notice in particular that the left hand side of condition \eqref{eq:main}, i.e., the SNR, is mainly a function of the desired controller $U_*^{\leq v}$ and the closed loop performance $\wt$ that it achieves, whereas the right hand side of \eqref{eq:main} is mainly a function of the structure of the optimal controller and $\Mtv$.  Thus we expect controllers with sparse and evenly distributed actuation, i.e., controllers that minimize the SNR threshold $(\gtv - \btv)^{-1}$, that act quickly and aggressively to achieve a good closed loop norm, i.e., controllers that maximize the SNR \eqref{eq:groupSNR}, to be recovered by the RFD procedure.

\section{Case Study}
\label{sec:case}
%

The following case study illustrates the concepts introduced in the previous section on a concrete system that satisfies our sufficient conditions.

\begin{figure}[h!]
\centering
\includegraphics[width = .35\textwidth]{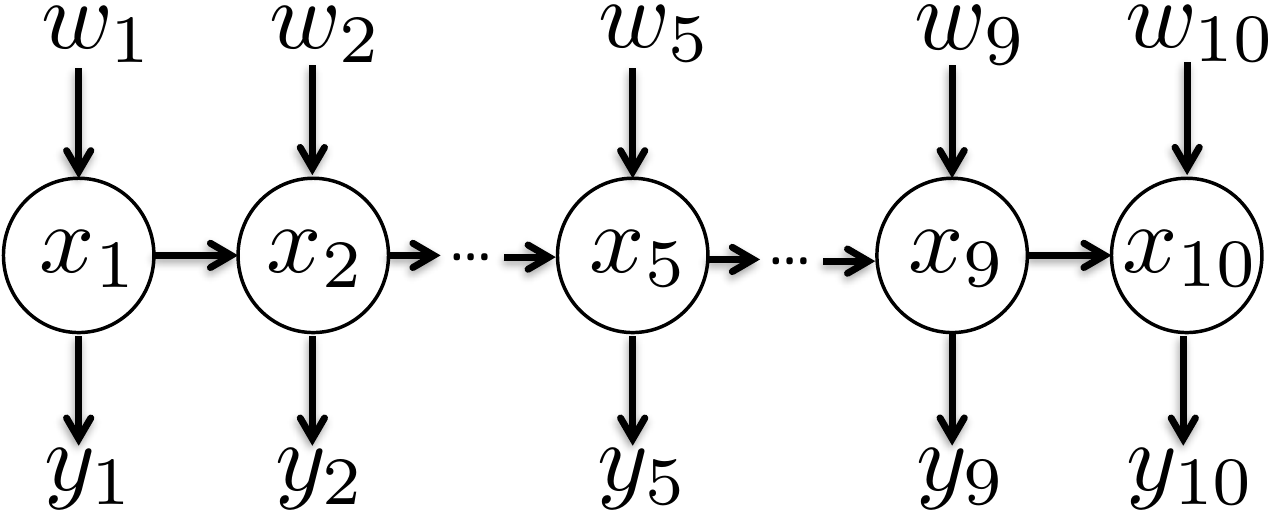}
\caption{A diagram of the Stable Unidirectional Chain System case study.}
\label{fig:chain}
\end{figure}

We consider a $\Htwo$ RFD optimization with control cost $\rho_u = .1$, and the remaining generalized plant \eqref{eq:genplant} state space parameters set as $B_2=C_1=I_{10}$, $A = \frac{1}{2}I_{10} + \frac{1}{2}Z_{10}$, and \edit{$B_1 = 1.1(E_{11} + E_{55}) + .7 E_{99} + .1 I_{10}$}.  This system is illustrated in Figure \ref{fig:chain}.  
This simple example is chosen in order to allow a direct computation of various bounds and parameters, and to easily interpret the propagation of inputs and disturbances.

We consider the task of recovering the optimal actuation schemes that use either 2 actuators or 3 actuators.    In particular, we take the desired controller $U_s$, for $s=2$ and $s=3$, to be 
\begin{equation}
\begin{array}{rl}
U_s := \argmin{U \in \RHinf}&\Htwonorm{Y-\M(U)}^2 + .1\Htwonorm{U}^2 \\
\text{s.t.} & \text{\edit{$U$ has at most $s$ nonzero rows}},
\end{array}
\label{eq:chainUs}
\end{equation}
with open loop state response $Y$ and map $\M$ as defined in Example \ref{ex:actuators}.  We solve this optimization problem by enumerating all possible actuation schemes, and we find that the optimal actuation scheme for $s=2$ is given by actuators at nodes $1$ and $5$, and for $s=3$ by actuators at nodes $1$, $5$ and $9$.

\edit{We emphasize that the goal of this case study is to illustrate the concepts introduced in the previous section, and to help the reader understand how the various parameters affect the recovery conditions -- in practice, $\m_*$ and $\Aa_*$ are not available.  Further, we note that the case study presented is, as far as we are aware, the first example in the literature of a system for which convex optimization provably identifies an optimal actuation architecture.}

\begin{figure}
 \begin{subfigure}{.3\columnwidth}
\includegraphics[width = \columnwidth]{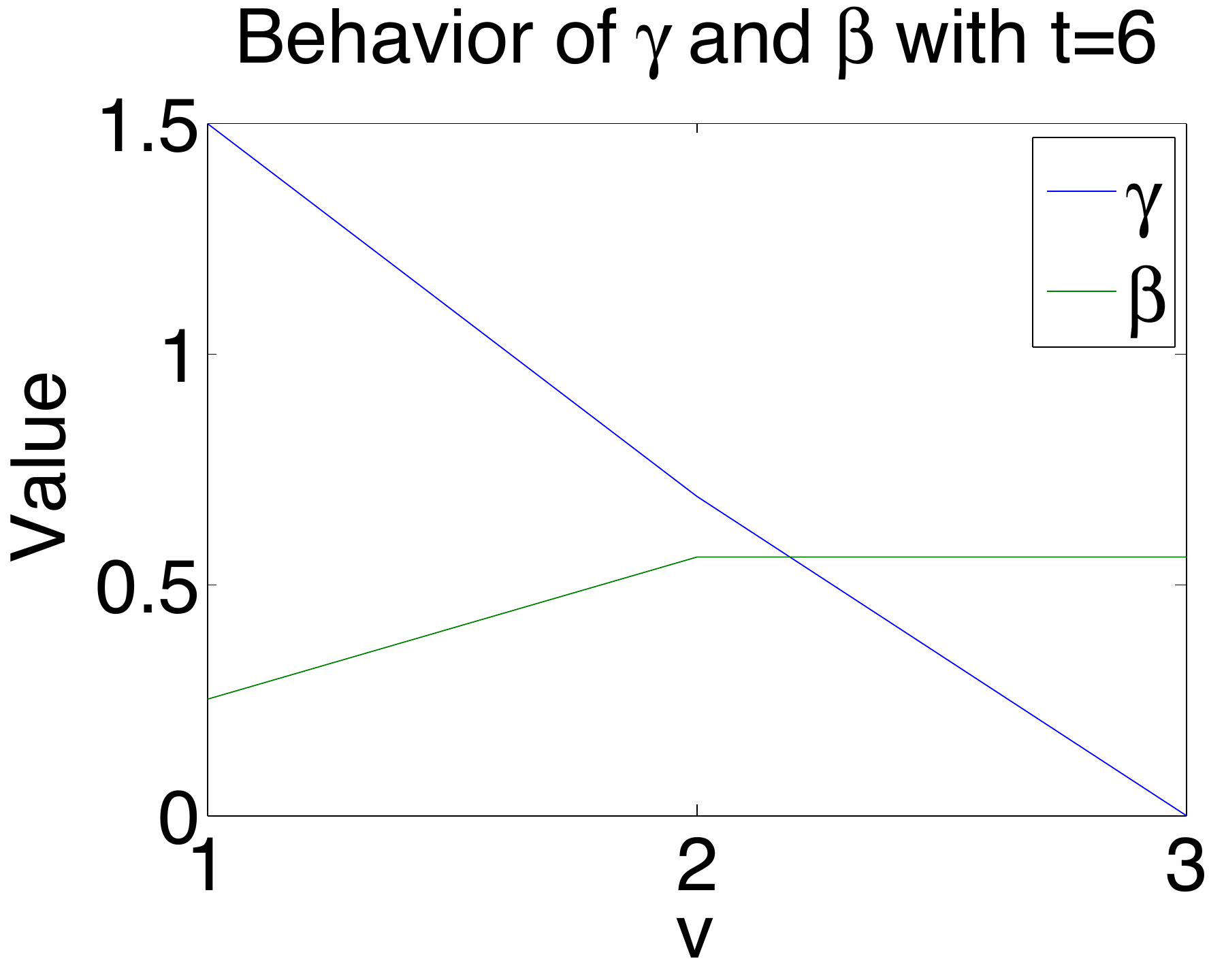}%
\caption{ Fixed horizon $t=6$.}%
\label{fig:GamBeta}%
\end{subfigure} \hfill %
\begin{subfigure}{.3\columnwidth}
\includegraphics[width = \columnwidth]{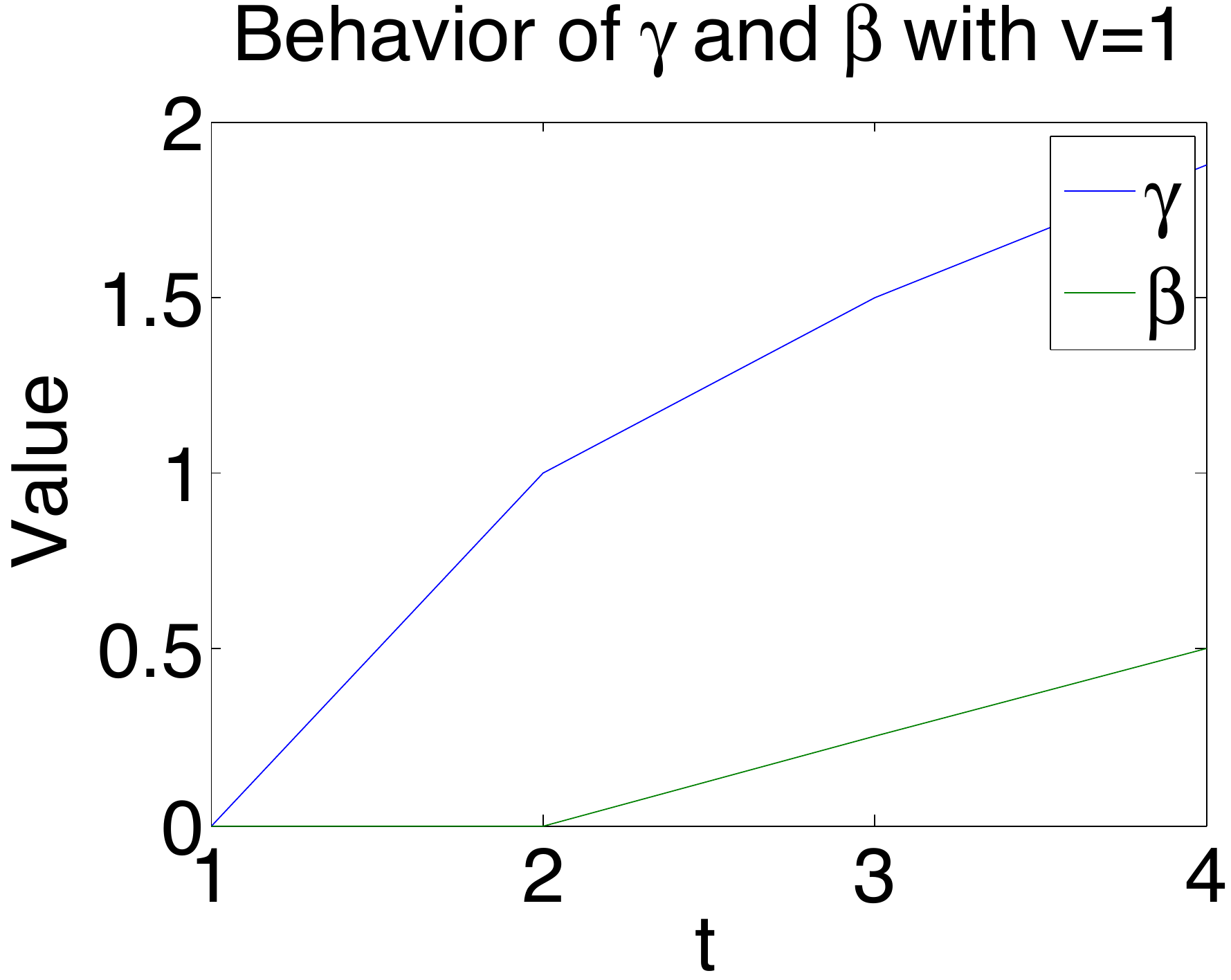}%
\caption{Fixed controller order $v=1$.}%
\label{fig:fixedV}%
\end{subfigure}\hfill%
\begin{subfigure}{.3\columnwidth}
\includegraphics[width = \columnwidth]{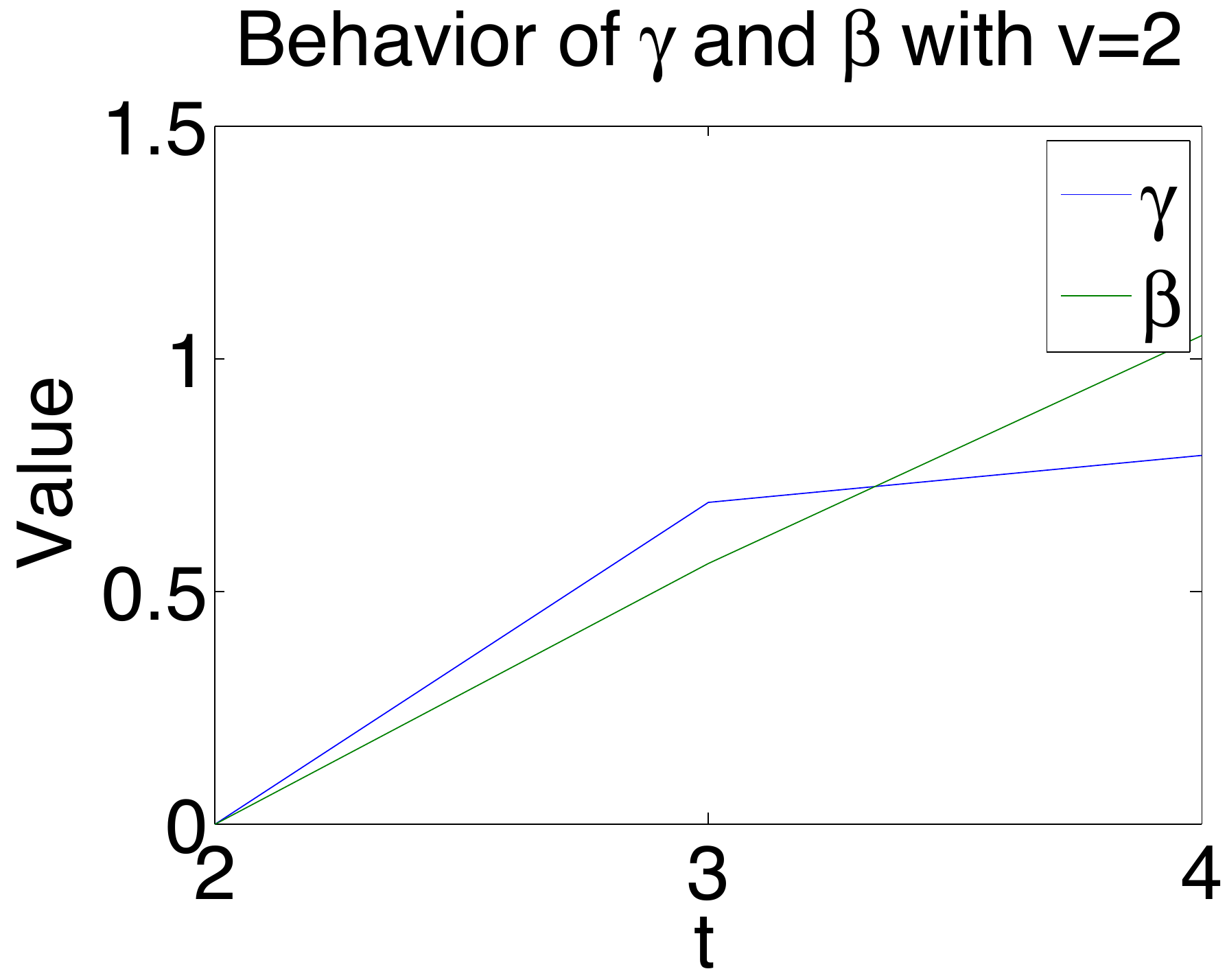}%
\caption{Fixed controller order $v=2$.}%
\label{fig:fixedV2}%
\end{subfigure} \hfill%
\caption{Behavior of identifiability parameters $\gtv$ and $\btv$.}
\end{figure}

With these optimal actuation schemes at our disposal, we vary the parameters $t$ and $v$ to investigate how our recovery conditions are affected.  As per the discussion in $\S\ref{sec:parameters}$, we set $\rho=0$.  We also show that for appropriate fixed controller order $v$ and horizon $t$, increasing $\lambda$ shifts the identified architecture from actuators at nodes 1, 5 and 9 to actuators at nodes 1 and 5.  This is a very desirable property from an architecture design perspective, as it suggests that increasing the regularization weight $\lambda$ causes the identified architecture to move to a simpler, but still optimal, actuator configuration.  \edit{As predicted by Corollary \ref{coro:success}, each optimal architecture is identified for a range of regularization weights $\lambda$.  Further, as predicted by Theorem \ref{thm:main}, the identified architectures are those for which the optimal control law achieves a small closed loop norm.}

We begin by examining how the lower bound parameter $\gtv$ and the {maximum gain} $\btv$ are affected as we vary the horizon $t$ and the controller order $v$.  In particular, for actuation sparsity $s=2$, we compute $\gtv$ and $\btv$ for (i) $t=6$ and $v\in\left\{1,2, 3\right\}$ (shown in Figure \ref{fig:GamBeta}),  (ii) for $v=1$ and $t\in\left\{1,2,3,4\right\}$ (shown in Figure \ref{fig:fixedV}), and (iii) for $v=2$ and $t\in\left\{2,3,4\right\}$ (shown in Figure \ref{fig:fixedV2}).  As expected, there is a decrease in the lower bound parameter $\gtv$ and an increase in the maximum gain $\btv$ as $v$ increases, while both $\gtv$ and $\btv$ are non-decreasing for a fixed controller order $v$ and increasing horizon $t$ as long as $t \leq \mix$. For this problem, the mixing time $\mix=5$.  Further we see that $\gtv$ begins to decrease for horizons $t> 6$ when $v=2$.


\begin{table}[h]
\centering
\begin{tabular}{c|cccccc}
\hline
 $t$ &  $\frac{1}{\gtv-\btv}$ & $\text{SNR}_1$ & $\text{SNR}_5$ & $\lambda$ &  $\Gdnorm{\Delta}$ & Bound \\ \hline
 2 & 1 & 1.27 & 1.27 & .8 & .73 & .89  \\ \hline
3&  .8 & .87 & .88 & 1.46 & .91 & 1.03 \\ \hline
 4 & .727 & .732 & .735 & 2.01 & 1.00 & 1.12  \\ \hline
5 & .7 & .67 & .68 & 2.45 & 1.05 & 1.16  \\ \hline
\end{tabular}
\caption{Summary of relevant values for the controller $U_2$ with actuators at nodes 1 and 5.}
\label{tab:chain2sparse}
\end{table}
\squeezeup
The conditions of Theorem \ref{thm:main} are satisfied for $v=1$ and several values of $t$.  For example, if we select $t=4$, $v=1$, $\rho = 0$ and $U_*$ as defined in \eqref{eq:chainUs}, we can compute $(\gamma^{\leq 4,1}-\beta^{\leq 4,1})^{-1}=.7273$, $\SNR{(U^{\leq 1}_*)_1} = .7324$, and $\SNR{(U^{\leq 1}_*)_5} = .7353$, thus satisfying condition \eqref{eq:main} for each of the two actuators.  Further, selecting $\lambda=2.0119 \in \Lambda$, and using this value for $\lambda$ in the truncated RFD optimization \eqref{eq:Topt} recovers a solution with non-zero first and fifth rows.


Perhaps surprisingly, similar positive recovery results can be verified for all $2\leq t \leq \mix$ -- the relevant values are summarized in Table \ref{tab:chain2sparse}.  In this table, $\text{SNR}_i$ corresponds to the SNR achieved by the controller component corresponding to the actuator at node $i$. Further, $\Delta := {\hat{U} - U^{\leq v}_*}$ is the approximation error between the architect parameter $\hat{U}$ and underlying parameter $U^{\leq v}_*$, and the values in the ``Bound'' column are given by equation \eqref{eq:bounded} in the appendix giving upper bounds on $\Gdnorm{\Delta}$.  It is worth noting that for $t=5$, we do not satisfy the sufficient conditions of Theorem \ref{thm:main}, but nonetheless recover the correct actuation architecture.  


We now consider the case of actuation sparsity $s=3$. Much as in the $s=2$ case, we can verify that the conditions of Theorem \ref{thm:main} hold for $v=1$, and $2\leq t \leq \mix$, where the mixing time $\mix$ is still 5.  However, since the controller with 3 actuators is able to achieve a much better closed loop norm, the SNRs are significantly larger, while the SNR threshold $(\gtv-\btv)^{-1}$ does not change significantly. In particular, for the case of $t=5$, we have a threshold of $(\gamma^{\leq 5,1}-\beta^{\leq 5,1})^{-1}=.82$, and SNRs of $4.04$, $4.04$ and $2.67$ for the three actuator components.  

This is consistent with our original interpretation of the closed loop state response $\wt$ playing the role of measurement noise -- the better the performance of the controller, the easier it is to identify via RFD.  These experiments demonstrate that controllers with sparse and diffuse actuation schemes that achieve a small state response $\wt$ are easy to identify as the solutions to RFD optimization problems.  In summary, our analysis and case studies demonstrate that: (i) the parameter $\rho$ can be set to 0 in the RFD optimization problem \eqref{eq:rfd_opt}, even if the original model matching problem \eqref{eq:gen_mm} had non-zero control cost $\rho_u$; (ii) choosing small controller order $v$ and horizon $t$ can actually lead to a favorable threshold $(\gtv - \btv)^{-1}$ \eqref{eq:main}; (iii) actuation schemes that are more evenly distributed (so that they lead to large mixing times $\mix$ in \eqref{eq:tau}) are easier to identify; and (iv) controller components that maximize the RFD analog of a SNR \eqref{eq:groupSNR} are more likely to satisfy our recovery conditions. These consist of controllers that have a concentration of energy in their early impulse response elements, and that achieve a closed loop with small state response component.

\section{Future Work}
\label{sec:conclusion}

\noindent\textbf{A priori bounds on incoherence}: It is of great interest to derive \apriori \ bounds on the gain parameters $\atv$ \eqref{eq:alpha} and $\btv$ \eqref{eq:beta} in terms of the state-space parameters of the system and a lower bound on the mixing time $\mix$ \eqref{eq:tau}.  This would allow for broader classes of systems to be classified as being amenable to RFD. \edit{We are currently pursuing semidefinite relaxation based methods to obtain bounds on these parameters \cite{YM_CDC15}.}

\noindent\edit{
\textbf{Scalability}: The scalability of the RFD framework is limited by the underlying quadratic invariance based controller synthesis algorithms upon which it is built.  In order to allow the RFD framework, and distributed optimal control theory in general, to scale to large heterogeneous systems, the first author and co-authors have developed the localized optimal control framework (cf. \cite{WM_CDC15llqg} and references therein).  The algorithmic component of the RFD framework has already been ported \cite{WMD_CDC15rfd}; it is of interest to see if analogous recovery conditions can also be developed.}   
\bibliographystyle{../IEEEtran}
\bibliography{biblio/comms,biblio/regularization_abr,biblio/matni_abr,biblio/sys_id_abr,biblio/decentralized_abr} 

\begin{thebibliography}{10}
\providecommand{\url}[1]{#1}
\csname url@rmstyle\endcsname
\providecommand{\newblock}{\relax}
\providecommand{\bibinfo}[2]{#2}
\providecommand\BIBentrySTDinterwordspacing{\spaceskip=0pt\relax}
\providecommand\BIBentryALTinterwordstretchfactor{4}
\providecommand\BIBentryALTinterwordspacing{\spaceskip=\fontdimen2\font plus
\BIBentryALTinterwordstretchfactor\fontdimen3\font minus
  \fontdimen4\font\relax}
\providecommand\BIBforeignlanguage[2]{{%
\expandafter\ifx\csname l@#1\endcsname\relax
\typeout{** WARNING: IEEEtran.bst: No hyphenation pattern has been}%
\typeout{** loaded for the language `#1'. Using the pattern for}%
\typeout{** the default language instead.}%
\else
\language=\csname l@#1\endcsname
\fi
#2}}

\bibitem{MC_CDC14}
N.~Matni and V.~Chandrasekaran, ``Regularization for design,'' \emph{CoRR},
  vol. arXiv:1404.1972, 2014.

\bibitem{Bon91}
F.~Bonsall, ``A general atomic decomposition theorem and banach's closed range
  theorem,'' \emph{The Quarterly Journal of Mathematics}, vol.~42, no.~1, pp.
  9--14, 1991.

\bibitem{CRT06}
E.~Candes, J.~Romberg, and T.~Tao, ``Robust uncertainty principles: exact
  signal reconstruction from highly incomplete frequency info.'' \emph{Info.
  Theory, IEEE Trans. on}, vol.~52, no.~2, pp. 489--509, Feb. 2006.

\bibitem{Donoho06}
D.~L. Donoho, ``Compressed sensing,'' \emph{IEEE Trans. Inform. Theory},
  vol.~52, pp. 1289--1306, 2006.

\bibitem{fazelThesis}
M.~Fazel, ``Matrix rank minimization with applications,'' Ph.D. dissertation,
  PhD thesis, Stanford University, 2002.

\bibitem{RFP10}
B.~Recht, M.~Fazel, and P.~A. Parrilo, ``Guaranteed minimum-rank solutions of
  linear matrix equations via nuclear norm minimization,'' \emph{SIAM Rev.},
  vol.~52, no.~3, pp. 471--501, Aug. 2010.

\bibitem{CRPW12}
V.~Chandrasekaran, B.~Recht, P.~Parrilo, and A.~Willsky,
  ``\BIBforeignlanguage{English}{The convex geometry of linear inverse
  problems},'' \emph{\BIBforeignlanguage{English}{Foundations of Computational
  Mathematics}}, vol.~12, pp. 805--849, 2012.

\bibitem{SBTR12}
P.~Shah, B.~N. Bhaskar, G.~Tang, and B.~Recht, ``Linear system identification
  via atomic norm regularization,'' in \emph{Decision and Control (CDC), 2012
  IEEE 51st Annual Conf. on}, 2012, pp. 6265--6270.

\bibitem{Fazel}
M.~Fazel, H.~Hindi, and S.~P. Boyd, ``A rank minimization heuristic with
  application to minimum order system approximation,'' in \emph{American
  Control Conference. Proc. of 2001}, vol.~6.\hskip 1em plus 0.5em minus
  0.4em\relax IEEE, 2001, pp. 4734--4739.

\bibitem{MR_CDC14}
N.~Matni and A.~Rantzer, ``Low-rank and low-order decompositions for local
  system identification,'' \emph{CoRR}, vol. arXiv:1403.7175, 2014.

\bibitem{LjungNewOld}
L.~Ljung, ``Some classical and some new ideas for identification of linear
  systems,'' \emph{Journal of Control, Automation and Electrical Systems},
  vol.~24, no. 1-2, pp. 3--10, 2013.

\bibitem{LFJ13}
F.~Lin, M.~Fardad, and M.~R. Jovanovic, ``Design of optimal sparse feedback
  gains via the alternating direction method of multipliers,'' \emph{Automatic
  Control, IEEE Trans. on}, vol.~58, no.~9, pp. 2426--2431, 2013.

\bibitem{JRM14}
V.~Jonsson, A.~Rantzer, and R.~M. Murray, ``A scalable formulation for
  engineering combination therapies for evolutionary dynamics of disease,'' in
  \emph{The IEEE American Control Conf. (ACC), 2014.}, 2014.

\bibitem{FLJ13}
M.~Fardad, F.~Lin, and M.~R. Jovanovi{\'c}, ``Design of optimal sparse
  interconnection graphs for synchronization of oscillator networks,''
  \emph{arXiv preprint arXiv:1302.0449}, 2013.

\bibitem{DJL14}
N.~Dhingra, M.~R. Jovanovic, and Z.-Q. Luo, ``An admm algorithm for optimal
  sensor and actuator selection,'' in \emph{IEEE Conf. on Decision and Control
  (CDC)}, 2014.

\bibitem{M_CDC13_codesign}
N.~Matni, ``Communication delay co-design in $\mathcal{H}_2$ decentralized
  control using atomic norm minimization,'' in \emph{Decision and Control
  (CDC), 2013 IEEE 52nd Annual Conf. on}, Dec 2013, pp. 6522--6529.

\bibitem{M_TCNS14}
------, ``Communication delay co-design in $\mathcal{H_{\mathrm{2}}}$
  distributed control using atomic norm minimization,'' \emph{IEEE Transactions
  on Networked Control Systems, Submitted to the}, vol. arXiv:1404.4911, 2015.

\bibitem{PKP15}
S.~Pequito, S.~Kar, and G.~J. Pappas, ``Minimum cost constrained input-output
  and control configuration co-design problem: A structural systems approach,''
  \emph{arXiv preprint arXiv:1503.02764}, 2015.

\bibitem{RL06}
M.~Rotkowitz and S.~Lall, ``A characterization of convex problems in
  decentralized control,'' \emph{Automatic Control, IEEE Trans. on}, vol.~51,
  no.~2, pp. 274--286, 2006.

\bibitem{RCL10}
M.~Rotkowitz, R.~Cogill, and S.~Lall, ``Convexity of optimal control over
  networks with delays and arbitrary topology,'' \emph{Int. J. Syst., Control
  Commun.}, vol.~2, no. 1/2/3, pp. 30--54, Jan. 2010.

\bibitem{LL11_QI}
L.~Lessard and S.~Lall, ``Quadratic invariance is necessary and sufficient for
  convexity,'' in \emph{American Control Conf. 2011}, 2011, pp. 5360--5362.

\bibitem{ZDG96}
K.~Zhou, J.~C. Doyle, K.~Glover, \emph{et~al.}, \emph{Robust and optimal
  control}.\hskip 1em plus 0.5em minus 0.4em\relax Prentice Hall New Jersey,
  1996, vol.~40.

\bibitem{LL11}
L.~Lessard and S.~Lall, ``A state-space solution to the two-player
  decentralized optimal control problem,'' in \emph{49th Annual Allerton Conf.
  on Comm., Control, and Computing}.\hskip 1em plus 0.5em minus 0.4em\relax
  IEEE, 2011, pp. 1559--1564.

\bibitem{S13}
C.~W. Scherer, ``Structured $\mathcal{H}_\infty$-optimal control for nested
  interconnections: A state-space solution,'' \emph{arXiv preprint
  arXiv:1305.1746}, 2013.

\bibitem{SP10}
P.~Shah and P.~A. Parrilo, ``$\mathcal{H}_2$-optimal decentralized control over
  posets: A state space solution for state-feedback,'' in \emph{Decision and
  Control (CDC), 2010 49th IEEE Conf. on}.\hskip 1em plus 0.5em minus
  0.4em\relax IEEE, 2010, pp. 6722--6727.

\bibitem{LL12}
A.~Lamperski and L.~Lessard, ``Optimal state-feedback control under sparsity
  and delay constraints,'' in \emph{3rd IFAC Workshop on Distributed Estimation
  and Control in Networked Systems}, 2012, pp. 204--209.

\bibitem{LD14}
A.~Lamperski and J.~C. Doyle, ``The $\mathcal{H}_2$ control problem for
  decentralized systems with delays,'' \emph{CoRR}, vol. abs/1312.7724, 2013.

\bibitem{M_CDC14_dhinf}
N.~Matni, ``Distributed control subject to delays satisfying an
  $\mathcal{H}_\infty$ norm bound,'' in \emph{Decision and Control (CDC), 2014
  IEEE Annual Conf. on}, Dec 2014, pp. 4006--4013.

\bibitem{HC72}
Y.-C. Ho and K.-C. Chu, ``Team decision theory and information structures in
  optimal control problems--part i,'' \emph{Automatic Control, IEEE Trans. on},
  vol.~17, no.~1, pp. 15--22, 1972.

\bibitem{DahlehL1}
M.~A. Dahleh and I.~J. Diaz-Bobillo, \emph{Control of uncertain systems: a
  linear programming approach}.\hskip 1em plus 0.5em minus 0.4em\relax
  Prentice-Hall, Inc., 1994.

\bibitem{SM12}
S.~Sabau and N.~C. Martins, ``Stabilizability and norm-optimal control design
  subject to sparsity constraints,'' \emph{arXiv:1209.1123}, 2012.

\bibitem{MMRY12}
A.~Mahajan, N.~Martins, M.~Rotkowitz, and S.~Yuksel, ``Information structures
  in optimal decentralized control,'' in \emph{Decision and Control (CDC), 2012
  IEEE 51st Annual Conf. on}, 2012, pp. 1291--1306.

\bibitem{Ridge}
A.~E. Hoerl and R.~W. Kennard, ``Ridge regression: Biased estimation for
  nonorthogonal problems,'' \emph{Technometrics}, vol.~12, pp. 55--67, 1970.

\bibitem{GroupLasso}
M.~Yuan and Y.~Lin, ``Model selection and estimation in regression with grouped
  variables,'' \emph{Journal of the Royal Statistical Society: Series B
  (Statistical Methodology)}, vol.~68, no.~1, pp. 49--67, 2006.

\bibitem{ElasticNet}
H.~Zou and T.~Hastie, ``Regularization and variable selection via the elastic
  net,'' \emph{Journal of the Royal Statistical Society: Series B (Statistical
  Methodology)}, vol.~67, no.~2, pp. 301--320, 2005.

\bibitem{OJV11}
G.~Obozinski, L.~Jacob, and J.-P. Vert, ``Group lasso with overlaps: the latent
  group lasso approach,'' \emph{arXiv preprint arXiv:1110.0413}, 2011.

\bibitem{RWY10}
G.~Raskutti, M.~J. Wainwright, and B.~Yu, ``Restricted eigenvalue properties
  for correlated gaussian designs,'' \emph{The Journal of Machine Learning
  Research}, vol.~11, pp. 2241--2259, 2010.

\bibitem{CT05}
E.~J. Candes and T.~Tao, ``Decoding by linear programming,'' \emph{Info.
  Theory, IEEE Trans. on}, vol.~51, no.~12, pp. 4203--4215, 2005.

\bibitem{YM_CDC15}
S.~You and N.~Matni, ``A convex approach to sparse h-infinity analysis and
  synthesis,'' in \emph{Decision and Control (CDC), 2015 IEEE Annual Conf. on,
  Submitted.}, 2015.

\bibitem{WM_CDC15llqg}
Y.-S. Wang, N.~Matni, and J.~C. Doyle, ``Localized lqg optimal control for
  large-scale systems,'' in \emph{Decision and Control (CDC), 2015 IEEE Annual
  Conf. on, Submitted.}, 2015.

\bibitem{WMD_CDC15rfd}
------, ``Localized lqr control with actuator regularization,'' in
  \emph{Decision and Control (CDC), 2015 IEEE Annual Conf. on, Submitted.},
  2015.

\end{thebibliography}


\begin{thebibliography}{10}
\providecommand{\url}[1]{#1}
\csname url@rmstyle\endcsname
\providecommand{\newblock}{\relax}
\providecommand{\bibinfo}[2]{#2}
\providecommand\BIBentrySTDinterwordspacing{\spaceskip=0pt\relax}
\providecommand\BIBentryALTinterwordstretchfactor{4}
\providecommand\BIBentryALTinterwordspacing{\spaceskip=\fontdimen2\font plus
\BIBentryALTinterwordstretchfactor\fontdimen3\font minus
  \fontdimen4\font\relax}
\providecommand\BIBforeignlanguage[2]{{%
\expandafter\ifx\csname l@#1\endcsname\relax
\typeout{** WARNING: IEEEtran.bst: No hyphenation pattern has been}%
\typeout{** loaded for the language `#1'. Using the pattern for}%
\typeout{** the default language instead.}%
\else
\language=\csname l@#1\endcsname
\fi
#2}}

\bibitem{MC_CDC14}
\BIBentryALTinterwordspacing
N.~Matni and V.~Chandrasekaran, ``Regularization for design,'' \emph{CoRR},
  vol. arXiv:1404.1972, 2014. [Online]. Available:
  \url{http://arxiv.org/abs/1404.1972}
\BIBentrySTDinterwordspacing

\bibitem{Bon91}
F.~Bonsall, ``A general atomic decomposition theorem and banach's closed range
  theorem,'' \emph{The Quarterly Journal of Mathematics}, vol.~42, no.~1, pp.
  9--14, 1991.

\bibitem{Pis86}
G.~Pisier, \emph{Probabilistic methods in the geometry of Banach spaces}.\hskip
  1em plus 0.5em minus 0.4em\relax Springer, 1986.

\bibitem{Jon92}
L.~K. Jones, ``A simple lemma on greedy approximation in hilbert space and
  convergence rates for projection pursuit regression and neural network
  training,'' \emph{The annals of Statistics}, pp. 608--613, 1992.

\bibitem{Bar93}
A.~R. Barron, ``Universal approximation bounds for superpositions of a
  sigmoidal function,'' \emph{Information Theory, IEEE Transactions on},
  vol.~39, no.~3, pp. 930--945, 1993.

\bibitem{DeVT96}
R.~A. DeVore and V.~N. Temlyakov, ``Some remarks on greedy algorithms,''
  \emph{Advances in computational Mathematics}, vol.~5, no.~1, pp. 173--187,
  1996.

\bibitem{CDS98}
S.~S. Chen, D.~L. Donoho, and M.~A. Saunders, ``Atomic decomposition by basis
  pursuit,'' \emph{SIAM journal on scientific computing}, vol.~20, no.~1, pp.
  33--61, 1998.

\bibitem{CRT06}
E.~Candes, J.~Romberg, and T.~Tao, ``Robust uncertainty principles: exact
  signal reconstruction from highly incomplete frequency information,''
  \emph{Information Theory, IEEE Transactions on}, vol.~52, no.~2, pp.
  489--509, Feb. 2006.

\bibitem{Donoho04}
D.~L. Donoho, ``For most large underdetermined systems of linear equations the
  minimal $\ell_1$-norm solution is also the sparsest solution,'' \emph{Comm.
  Pure Appl. Math}, vol.~59, pp. 797--829, 2004.

\bibitem{Donoho06}
------, ``Compressed sensing,'' \emph{IEEE Trans. Inform. Theory}, vol.~52, pp.
  1289--1306, 2006.

\bibitem{fazelThesis}
M.~Fazel, ``Matrix rank minimization with applications,'' Ph.D. dissertation,
  PhD thesis, Stanford University, 2002.

\bibitem{RFP10}
B.~Recht, M.~Fazel, and P.~A. Parrilo, ``Guaranteed minimum-rank solutions of
  linear matrix equations via nuclear norm minimization,'' \emph{SIAM Rev.},
  vol.~52, no.~3, pp. 471--501, Aug. 2010.

\bibitem{CR12}
E.~Cand\`{e}s and B.~Recht, ``Exact matrix completion via convex
  optimization,'' \emph{Commun. ACM}, vol.~55, no.~6, pp. 111--119, June 2012.

\bibitem{CRPW12}
V.~Chandrasekaran, B.~Recht, P.~Parrilo, and A.~Willsky,
  ``\BIBforeignlanguage{English}{The convex geometry of linear inverse
  problems},'' \emph{\BIBforeignlanguage{English}{Foundations of Computational
  Mathematics}}, vol.~12, pp. 805--849, 2012.

\bibitem{SBTR12}
P.~Shah, B.~N. Bhaskar, G.~Tang, and B.~Recht, ``Linear system identification
  via atomic norm regularization,'' in \emph{Decision and Control (CDC), 2012
  IEEE 51st Annual Conference on}, 2012, pp. 6265--6270.

\bibitem{Fazel}
M.~Fazel, H.~Hindi, and S.~P. Boyd, ``A rank minimization heuristic with
  application to minimum order system approximation,'' in \emph{American
  Control Conference, 2001. Proceedings of the 2001}, vol.~6.\hskip 1em plus
  0.5em minus 0.4em\relax IEEE, 2001, pp. 4734--4739.

\bibitem{MR_CDC14}
\BIBentryALTinterwordspacing
N.~Matni and A.~Rantzer, ``Low-rank and low-order decompositions for local
  system identification,'' \emph{CoRR}, vol. arXiv:1403.7175, 2014. [Online].
  Available: \url{http://arxiv.org/abs/1403.7175}
\BIBentrySTDinterwordspacing

\bibitem{LjungNewOld}
L.~Ljung, ``Some classical and some new ideas for identification of linear
  systems,'' \emph{Journal of Control, Automation and Electrical Systems},
  vol.~24, no. 1-2, pp. 3--10, 2013.

\bibitem{LFJ13}
F.~Lin, M.~Fardad, and M.~R. Jovanovic, ``Design of optimal sparse feedback
  gains via the alternating direction method of multipliers,'' \emph{Automatic
  Control, IEEE Transactions on}, vol.~58, no.~9, pp. 2426--2431, 2013.

\bibitem{JRM14}
V.~Jonsson, A.~Rantzer, and R.~M. Murray, ``A scalable formulation for
  engineering combination therapies for evolutionary dynamics of disease,'' in
  \emph{The IEEE American Control Conference (ACC), 2014.}, 2014.

\bibitem{DLFM12}
N.~Dhingra, F.~Lin, M.~Fardad, and M.~R. Jovanovic, ``On identifying sparse
  representations of consensus networks,'' in \emph{3rd IFAC Workshop on
  Distributed Estimation and Control in Networked Systems, Santa Barbara, CA},
  2012, pp. 305--310.

\bibitem{XB07}
L.~Xiao, S.~Boyd, and S.-J. Kim, ``Distributed average consensus with
  least-mean-square deviation,'' \emph{Journal of Parallel and Distributed
  Computing}, vol.~67, no.~1, pp. 33--46, 2007.

\bibitem{FLJ13}
M.~Fardad, F.~Lin, and M.~R. Jovanovi{\'c}, ``Design of optimal sparse
  interconnection graphs for synchronization of oscillator networks,''
  \emph{arXiv preprint arXiv:1302.0449}, 2013.

\bibitem{DJL14}
N.~Dhingra, M.~R. Jovanovic, and Z.-Q. Luo, ``An admm algorithm for optimal
  sensor and actuator selection,'' in \emph{IEEE Conference on Decision and
  Control (CDC), 2014, To appear in the proceedings of the}, 2014.

\bibitem{M_CDC13_codesign}
N.~Matni, ``Communication delay co-design in $\mathcal{H}_2$ decentralized
  control using atomic norm minimization,'' in \emph{Decision and Control
  (CDC), 2013 IEEE 52nd Annual Conference on}, Dec 2013, pp. 6522--6529.

\bibitem{M_TCNS14}
\BIBentryALTinterwordspacing
------, ``Communication delay co-design in $\mathcal{H_{\mathrm{2}}}$
  distributed control using atomic norm minimization,'' \emph{IEEE Transactions
  on Networked Control Systems, Submitted to the}, vol. arXiv:1404.4911, 2014.
  [Online]. Available: \url{http://arxiv.org/abs/1404.4911}
\BIBentrySTDinterwordspacing

\bibitem{PKA15}
S.~Pequito, S.~Kar, and A.~Aguiar, ``A framework for structural input/output
  and control configuration selection in large-scale systems,'' \emph{Automatic
  Control, IEEE Transactions on}, vol.~PP, no.~99, pp. 1--1, 2015.

\bibitem{PKP15}
S.~Pequito, S.~Kar, and G.~J. Pappas, ``Minimum cost constrained input-output
  and control configuration co-design problem: A structural systems approach,''
  \emph{arXiv preprint arXiv:1503.02764}, 2015.

\bibitem{RL06}
M.~Rotkowitz and S.~Lall, ``A characterization of convex problems in
  decentralized control,'' \emph{Automatic Control, IEEE Transactions on},
  vol.~51, no.~2, pp. 274--286, 2006.

\bibitem{RCL10}
\BIBentryALTinterwordspacing
M.~Rotkowitz, R.~Cogill, and S.~Lall, ``Convexity of optimal control over
  networks with delays and arbitrary topology,'' \emph{Int. J. Syst., Control
  Commun.}, vol.~2, no. 1/2/3, pp. 30--54, Jan. 2010. [Online]. Available:
  \url{http://dx.doi.org/10.1504/IJSCC.2010.031157}
\BIBentrySTDinterwordspacing

\bibitem{LL11_QI}
L.~Lessard and S.~Lall, ``Quadratic invariance is necessary and sufficient for
  convexity,'' in \emph{American Control Conference (ACC), 2011}.\hskip 1em
  plus 0.5em minus 0.4em\relax IEEE, 2011, pp. 5360--5362.

\bibitem{ZDG96}
K.~Zhou, J.~C. Doyle, K.~Glover, \emph{et~al.}, \emph{Robust and optimal
  control}.\hskip 1em plus 0.5em minus 0.4em\relax Prentice Hall New Jersey,
  1996, vol.~40.

\bibitem{LL11}
L.~Lessard and S.~Lall, ``A state-space solution to the two-player
  decentralized optimal control problem,'' in \emph{49th Annual Allerton
  Conference on Communication, Control, and Computing}.\hskip 1em plus 0.5em
  minus 0.4em\relax IEEE, 2011, pp. 1559--1564.

\bibitem{S13}
C.~W. Scherer, ``Structured $\mathcal{H}_\infty$-optimal control for nested
  interconnections: A state-space solution,'' \emph{arXiv preprint
  arXiv:1305.1746}, 2013.

\bibitem{SP10}
P.~Shah and P.~A. Parrilo, ``$\mathcal{H}_2$-optimal decentralized control over
  posets: A state space solution for state-feedback,'' in \emph{Decision and
  Control (CDC), 2010 49th IEEE Conference on}.\hskip 1em plus 0.5em minus
  0.4em\relax IEEE, 2010, pp. 6722--6727.

\bibitem{TP14}
T.~Tanaka and P.~A. Parrilo, ``Optimal output feedback architecture for
  triangular {LQG} problems,'' \emph{arXiv preprint arXiv:1403.4330}, 2014.

\bibitem{LL12}
A.~Lamperski and L.~Lessard, ``Optimal state-feedback control under sparsity
  and delay constraints,'' in \emph{3rd IFAC Workshop on Distributed Estimation
  and Control in Networked Systems}, 2012, pp. 204--209.

\bibitem{LD14}
\BIBentryALTinterwordspacing
A.~Lamperski and J.~C. Doyle, ``The $\mathcal{H}_2$ control problem for
  decentralized systems with delays,'' \emph{CoRR}, vol. abs/1312.7724, 2013.
  [Online]. Available: \url{http://arxiv.org/abs/1312.7724}
\BIBentrySTDinterwordspacing

\bibitem{LDXX}
------, ``Output feedback $\mathcal{H}_2$ model matching for decentralized
  systems with delays,'' in \emph{American Control Conference (ACC),
  2013}.\hskip 1em plus 0.5em minus 0.4em\relax IEEE, 2013, pp. 5778--5783.

\bibitem{M_CDC14_dhinf}
N.~Matni, ``Distributed control subject to delays satisfying an
  $\mathcal{H}_\infty$ norm bound,'' in \emph{Decision and Control (CDC), 2014
  IEEE 53rd Annual Conference on}, Dec 2014, pp. 4006--4013.

\bibitem{Ridge}
A.~E. Hoerl and R.~W. Kennard, ``Ridge regression: Biased estimation for
  nonorthogonal problems,'' \emph{Technometrics}, vol.~12, no.~1, pp. 55--67,
  1970.

\bibitem{GroupLasso}
M.~Yuan and Y.~Lin, ``Model selection and estimation in regression with grouped
  variables,'' \emph{Journal of the Royal Statistical Society: Series B
  (Statistical Methodology)}, vol.~68, no.~1, pp. 49--67, 2006.

\bibitem{ElasticNet}
H.~Zou and T.~Hastie, ``Regularization and variable selection via the elastic
  net,'' \emph{Journal of the Royal Statistical Society: Series B (Statistical
  Methodology)}, vol.~67, no.~2, pp. 301--320, 2005.

\bibitem{SM12}
S.~Sabau and N.~C. Martins, ``Stabilizability and norm-optimal control design
  subject to sparsity constraints,'' \emph{arXiv:1209.1123}, 2012.

\bibitem{MMRY12}
A.~Mahajan, N.~Martins, M.~Rotkowitz, and S.~Yuksel, ``Information structures
  in optimal decentralized control,'' in \emph{Decision and Control (CDC), 2012
  IEEE 51st Annual Conference on}, 2012, pp. 1291--1306.

\bibitem{HC72}
Y.-C. Ho and K.-C. Chu, ``Team decision theory and information structures in
  optimal control problems--part i,'' \emph{Automatic Control, IEEE
  Transactions on}, vol.~17, no.~1, pp. 15--22, 1972.

\bibitem{DahlehL1}
M.~A. Dahleh and I.~J. Diaz-Bobillo, \emph{Control of uncertain systems: a
  linear programming approach}.\hskip 1em plus 0.5em minus 0.4em\relax
  Prentice-Hall, Inc., 1994.

\bibitem{JOV09}
L.~Jacob, G.~Obozinski, and J.-P. Vert, ``Group lasso with overlap and graph
  lasso,'' in \emph{Proceedings of the 26th Annual International Conference on
  Machine Learning}, ser. ICML '09.\hskip 1em plus 0.5em minus 0.4em\relax New
  York, NY, USA: ACM, 2009, pp. 433--440.

\bibitem{OJV11}
G.~Obozinski, L.~Jacob, and J.-P. Vert, ``Group lasso with overlaps: the latent
  group lasso approach,'' \emph{arXiv preprint arXiv:1110.0413}, 2011.

\bibitem{RWY10}
G.~Raskutti, M.~J. Wainwright, and B.~Yu, ``Restricted eigenvalue properties
  for correlated gaussian designs,'' \emph{The Journal of Machine Learning
  Research}, vol.~11, pp. 2241--2259, 2010.

\bibitem{VB09}
S.~A. Van De~Geer, P.~B{\"u}hlmann, \emph{et~al.}, ``On the conditions used to
  prove oracle results for the lasso,'' \emph{Electronic Journal of
  Statistics}, vol.~3, pp. 1360--1392, 2009.

\bibitem{CT05}
E.~J. Candes and T.~Tao, ``Decoding by linear programming,'' \emph{Information
  Theory, IEEE Transactions on}, vol.~51, no.~12, pp. 4203--4215, 2005.

\bibitem{d2007direct}
A.~d'Aspremont, L.~El~Ghaoui, M.~I. Jordan, and G.~R. Lanckriet, ``A direct
  formulation for sparse pca using semidefinite programming,'' \emph{SIAM
  review}, vol.~49, no.~3, pp. 434--448, 2007.

\bibitem{YM_CDC15}
S.~You and N.~Matni, ``A convex approach to sparse h-infinity analysis and
  synthesis,'' in \emph{Decision and Control (CDC), 2015 IEEE 54th Annual
  Conference on, Submitted.}, 2015.

\bibitem{TL11}
T.~Tanaka and C.~Langbort, ``The bounded real lemma for internally positive
  systems and h-infinity structured static state feedback,'' \emph{IEEE
  transactions on automatic control}, vol.~56, no.~9, pp. 2218--2223, 2011.

\bibitem{R12}
A.~Rantzer, ``Distributed control of positive systems,'' \emph{arXiv preprint
  arXiv:1203.0047}, 2012.

\bibitem{DTF13}
K.~Dvijotham, E.~Todorov, and M.~Fazel, ``Convex structured controller
  design,'' \emph{arXiv preprint arXiv:1309.7731}, 2013.

\bibitem{WMYD_ACC14}
Y.-S. Wang, N.~Matni, S.~You, and J.~C. Doyle, ``Localized distributed state
  feedback control with communication delays,'' in \emph{American Control
  Conference (ACC), 2014}.\hskip 1em plus 0.5em minus 0.4em\relax IEEE, 2014,
  pp. 5748--5755.

\bibitem{WMD_CDC14}
\BIBentryALTinterwordspacing
Y.-S. Wang, N.~Matni, and J.~C. Doyle, ``Localized lqr optimal control,'' in
  \emph{Decision and Control (CDC), 2014 IEEE 53rd Annual Conference on, To
  appear.}, 2014. [Online]. Available: \url{http://arxiv.org/abs/1409.6404}
\BIBentrySTDinterwordspacing

\bibitem{WM_Allerton14}
Y.-S. Wang and N.~Matni, ``Localized distributed optimal control with output
  feedback and communication delays,'' in \emph{Communication, Control, and
  Computing, IEEE 52nd Annual Allerton Conference on}, 2014.

\end{thebibliography}
\vspace{-1.35mm}
\begin{appendix}
\section{Proof of Theorem \ref{thm:consistency}}
\begin{IEEEproof}
The proof of Theorem \ref{thm:consistency} centers around showing that under Assumption \ref{as:nu}, the {unique} solution to the architect optimization \eqref{eq:Toracle}
is also the {unique} solution of the original unconstrained optimization \eqref{eq:Topt}.  We emphasize that at no point during the RFD process do we assume knowledge of $\Gs$ or of the architect optimization problem \eqref{eq:Toracle}.  

The proof consists of two parts: we first show that if $\atv>0$, the architect optimization problem \eqref{eq:Toracle} has a unique optimal solution $\hat{U}$, and control its deviation from the underlying desired controller $U^{\leq v}_*$.  We then use $\hat{U}$ and its error bound to construct a strictly dual-feasible primal/dual pair for the {original} RFD optimization problem \eqref{eq:Topt}, showing that $\hat{U}$ is indeed its unique optimal solution as well.

\begin{prop}[Bounded Errors]
Fix a horizon $1 \leq t < \infty$, and a controller order $1 \leq v \leq t$.  Assume that $\atv$ as defined in \eqref{eq:alpha} is strictly positive, and let $\Delta := \hat{U} - {U_*^{\leq v}}$.  Then
\vspace{-2mm}
\begin{equation}
\begin{array}{l}
\Gdnorm{\Delta}\leq \frac{1}{\alpha}\left(\lambda + \Gdnorm{\adj{\Mtv_{\Gs}} (\err)}+ \rho \Gdnorm{{U_*^{\leq v}}}\right).
\end{array}
\label{eq:bounded}
\end{equation}
\squeezeup
\label{prop:bounded}
\end{prop}
\begin{IEEEproof}
It is clear that under the assumption that $\atv > 0$, the architect optimization problem \eqref{eq:Toracle} is strongly convex, and hence has a unique optimal solution $\hat{U}$. 
Letting $\Delta:=\hat{U}-{U^{\leq v}_*}$, and using the relation \eqref{eq:Ttruth}, the optimality conditions of the architect optimization problem \eqref{eq:Toracle} are then given by
\vspace{-2mm}
\begin{equation*}
\begin{array}{l}
\left(\adj{\Mtva} \Mtva + \rho I\right)(\Delta) - \adj{\Mtva}(\err)+ \rho {U^{\leq v}_*} + \lambda Z + \Lambda_{\Gs^\perp} \ni 0,
\label{eq:opt_cdt}
\end{array}
\end{equation*}
where $Z \in \partial \Gnorm{\hat{U}}$ satisfies $\Gdnorm{Z_{\Gs}}=1$, $\Gdnorm{Z_{\Gs^\perp}}\leq 1$, and $\Lambda_{\Gs^\perp} \in \Gs^\perp$ is the Lagrange multiplier corresponding to the architect constraint $U \in \Gs$.  Projecting \eqref{eq:opt_cdt} onto $\Gs$, and leveraging that $\Delta \in \Gs$, we then obtain
\begin{multline}
\left(\adj{\Mtva} \Mtva + \rho I\right)(\Delta) =  \left(\adj{\Mtva}(\err) - \rho {U^{\leq v}_*} - \lambda Z_{\Gs}\right).
\label{eq:pGs}
\end{multline}

 We then have the following chain of inequalities
\begin{equation*}
\begin{array}{l}
\atv\Gdnorm{\Delta} \leq  \Gdnorm{\left(\adj{\Mtva} \Mtva + \rho I\right)(\Delta)} 
\leq \lambda + \Gdnorm{\adj{\Mtva}(\err)} + \rho \Gdnorm{{U^{\leq v}_*}}
\end{array}
\end{equation*}
where the first inequality follows from \eqref{eq:alpha}, and the second from \eqref{eq:pGs} and the triangle inequality.  Rearranging terms yields the error bound \eqref{eq:bounded}.
\end{IEEEproof}

\subsubsection*{Strict dual feasibility} In order to construct a primal/dual feasible pair for optimization \eqref{eq:Topt} from $\hat{U}$, we first set $Z_{\Gs}$ to be a member of the sub differential $\partial \Gnorm{\cdot}$ evaluated at $\hat{U}$.  We now choose $Z_{\Gs^\perp}$ to be
\begin{equation}
\begin{array}{l}
Z_{\Gs^\perp} := \frac{\left( \adj{\Mtvap}(\err) -  \adj{\Mtvap} \Mtva(\Delta)\right)}{\lambda}
\end{array}
\label{eq:ZGsperp}
\end{equation}

In doing so, we guarantee that $(\hat{U}, Z)$ satisfy the optimality conditions of optimization \eqref{eq:Topt}.  What remains to be shown is the $Z_{\Gs^\perp}$ is an element of the sub-differential.  In order to do so, we show that under the assumptions of the theorem, $\Gdnorm{Z_{\Gs^\perp}} < 1$.  This guarantees that $Z$ is indeed in $\partial \Gnorm{\hat{U}}$, and that $\hat{U}_\Aa = 0$ for all $\Aa \notin \m_*$.  

To that end, notice that $\Gdnorm{Z_{\Gs^\perp}}$ can be upper bounded by
\begin{equation*}
\begin{array}{l}
\frac{\left(\Gdnorm{\adj{\Mtvap} \Mtva(\Delta)} + \Gdnorm{\adj{\Mtvap}(W^{\leq t}+\drv)} \right)}{\lambda} \\
\leq  \frac{1}{\lambda} \left(\btv\Gdnorm{\Delta} + \Gdnorm{\adj{\Mtvap}(W^{\leq t}+\drv)} \right) \\
 \leq \frac{1}{\lambda} \nu\left(\rho \Gdnorm{U^{\leq v}_*} + \Gdnorm{\adj{\Mtva}(W^{\leq t}+\drv)} \right)  \\
+ \nu + \frac{1}{\lambda} \Gdnorm{\adj{\Mtvap}(W^{\leq t}+\drv)}
 <  1,\end{array}
\end{equation*}
where the first inequality follows from applying the triangle inequality to \eqref{eq:ZGsperp}, the second from applying definition \eqref{eq:beta}, the third from applying the error bound \eqref{eq:bounded}, and the fourth from \eqref{eq:sufficient}.
Thus we have shown that under the assumptions of the Theorem, $\hat{U}$ is also the optimal solution of the original problem \eqref{eq:Topt} -- its uniqueness follows from the \emph{local} strong convexity of the cost function around $\hat{U}$.  
Finally, if for $\Aa \in \m_*$, we have that \eqref{eq:structure} holds, then $\hat{U}_\Aa \neq 0$.
\end{IEEEproof}
\begin{IEEEproof}[Proof of Lemma \ref{lem:param_bounds}]
The gain $\atv$ is bounded below by
\begin{equation*}
\begin{array}{l}
 \displaystyle \min_{\begin{array}{c}\Gdnorm{\Delta} = 1 \\ \Delta \in \Gs \end{array}} \max_{\Aa \in \m_*} \Htwonorm{\left(\adj{\Mtv_{\Aa}} \Mtv_\Aa + \rho I\right) \Delta_\Aa}  -\sum_{\Bb\neq \Aa \in \m_*} \Htwonorm{\adj{\Mtv_{\Aa}} \Mtv_\Bb \Delta_\Bb} \\ \indent
\geq \rho +\displaystyle \min_{\B \subseteq \m_*, |\B|\geq 1}\max_{\Aa\in \B} 
 \sm{\adj{\Mtv_{\Aa}} \Mtv_\Aa} - \sum_{\Bb\neq \Aa \in \B} \sM{\adj{\Mtv_{\Aa}} \Mtv_\Bb},
 \end{array}
\end{equation*}
where the inequalities follow from the fact that $\Gdnorm{\Delta} = 1$ implies that there exists $\Aa \in \m_*$ such that $\Htwonorm{\Delta_\Aa} = 1$, and the definition of the respective norms.  The derivation of the bound on $\btv$ is similar, and hence omitted.
\end{IEEEproof}
\end{appendix}

\end{document}